\documentclass[12pt,leqno]{article}

\usepackage{amsmath,amssymb,latexsym,theorem,epsfig}
\usepackage[all]{xy}
\CompileMatrices

\newcommand{\cp}[1]{{\mathbb P}^{#1}}
\newcommand{\op}[1]{\operatorname{#1}}

\newcommand{\FM}{\boldsymbol{F}{\boldsymbol{M}}}
\newcommand{\fm}{\boldsymbol{f}\boldsymbol{m}}
\newcommand{\T}{\boldsymbol{T}}

\newcommand{\num}[2]{$({\#}\text{{\bfseries #1}}#2)$}

\newcommand{\cV}{{\mathcal V}}

\newcommand{\heck}{\boldsymbol{H}\boldsymbol{e}\boldsymbol{c}
\boldsymbol{k}\boldsymbol{e}}
\newcommand{\hup}[2]{\heck^{+}_{(#1)}(#2)}
\newcommand{\hdown}[2]{\heck^{-}_{(#1)}(#2)}
\newcommand{\homsh}[3]{{\mathcal H}om_{#1}(#2,#3)}
\newcommand{\extsh}[4]{{\mathcal E}xt_{#2}^{#1}(#3,#4)}
\newcommand{\p}[1]{{\mathbb P}(#1)}
\newcommand{\bl}[2]{\operatorname{Bl}_{#1}#2}

\newtheorem{theo}{Theorem}[section]
\newtheorem{lem}[theo]{Lemma}
\newtheorem{cor}[theo]{Corollary}
\newtheorem{prop}[theo]{Proposition}
\newtheorem{defi}[theo]{Definition}
\newtheorem{claim}[theo]{Claim}

{\theorembodyfont{\rmfamily} \newtheorem{rem}[theo]{Remark}}
{\theorembodyfont{\rmfamily} }

\def\punkt{\refstepcounter{subsubsection} \medskip
           \noindent{\bf \thesubsubsection.\ }}

\numberwithin{equation}{section}
\newcounter{rom}

\newcommand{\Appendix}[1]{%
  \refstepcounter{section}%
  \addcontentsline{toc}{section}%
    {\bfseries\appendixname~\thesection\ #1}%
    {\medskip\noindent \Large\bfseries\appendixname\ \thesection\ #1}%
\sectionmark{#1}\smallskip\noindent
\renewcommand{\theequation}{\thesection.\arabic{equation}}
}

\begin{document}

\title{  \vspace*{-5em} 
\  \hfill {\normalsize CERN-TH/2000-203, UPR-894T, RU-00-5B} \\[0.5em]
\ \\ 
{\sf\LARGE Standard-model bundles}}
\author{Ron Donagi$^1$, Burt A.~Ovrut$^2$, Tony Pantev$^1$ 
      and Daniel Waldram$^{3}$ \\ [0.5em]
      \ \\
   {\normalsize $^1$Department of Mathematics, 
      University of Pennsylvania} \\[-0.2em]
      {\normalsize Philadelphia, PA 19104--6395, USA}\\
   {\normalsize $^2$Department of Physics, 
      University of Pennsylvania} \\[-0.2em]
      {\normalsize Philadelphia, PA 19104--6396, USA}\\
    {\normalsize $^3$Theory Division, CERN CH-1211, Geneva 23, 
Switzerland, and }\\[-0.2em]
{\normalsize Department of Physics, The Rockfeller University} \\[-0.2em]
{\normalsize New York, NY 10021}
}
\date{}
\maketitle

\begin{abstract}
We describe a family of genus one fibered Calabi-Yau threefolds with
fundamental group ${\mathbb Z}/2$. On each Calabi-Yau $Z$ in the
family we exhibit a
positive dimensional family of Mumford stable bundles whose symmetry
group is the Standard Model group $SU(3)\times SU(2)\times U(1)$ and
which have $c_{3} = 6$. We also show that for each bundle $V$
in our family, $c_{2}(Z) - c_{2}(V)$ is the class of an effective
curve on $Z$. These conditions ensure that $Z$ and $V$ can be used for a
phenomenologically relevant compactification of Heterotic M-theory.

\smallskip

\noindent
{\bf MSC 2000:} 14D21, 14J32
\end{abstract}

\section{Introduction} \label{s-intro}

In this paper we construct a particular class of bundles with
constrained Chern classes on certain non-simply connected Calabi-Yau
threefolds. These bundles are instrumental in deriving the Standard
Model of particle physics in the context of the Heterotic M-theory
\cite{usnew}. Bundles of this type 
have been the subject of active research for quite some
time \cite{tian-yau},
\cite{kachru},  \cite{pokorski-ross}, \cite{ACK}, \cite{usold},
\cite{richardt}. In contrast with the classical  
constructions \cite{tian-yau},
\cite{kachru},  \cite{pokorski-ross}, where the bundles obtained are
associated with 
the tangent bundle of the Calabi-Yau manifold and tend to be rigid, 
our examples are
independent of the geometry of the tangent bundle and vary in
families. In particular we construct infinitely many positive
dimensional families of bundles which are suitable for
phenomenologically relevant compactifications of Heterotic
M-theory. Our construction takes place entirely within the realm of 
algebraic geometry. 
The physical implications of our results are discussed in the
companion paper \cite{usnew}, which also contains a summary of the
construction written for physicists. In the remainder of this
introduction we give a brief overview of the physical motivation for
our work, followed by an outline of the actual geometric construction.

The search for exotic principal bundles on Calabi-Yau threefolds is
motivated by string theory. To compactify the $E_{8}\times
E_{8}$-Heterotic string to four dimensions  one
prescribes:
\begin{itemize}
\item a Calabi-Yau 3-fold $Z$; 
\item a Ricci flat K\"{a}hler metric on $Z$ with a K\"{a}hler form
$\omega$; 
\item an $\omega$-instanton ${\mathcal E} \to Z$ with a structure group
$E_{8}\times E_{8}$.
\end{itemize}

\medskip

The Hermit-Einstein connection on ${\mathcal E}$ is a vacuum of 
the Heterotic string theory. The moduli space of
${\mathcal E}$'s is a subspace of of the moduli space of vacua for the
Heterotic string. In view of the Uhlenbeck-Yau theorem
\cite{uhlenbeck-yau} every such ${\mathcal
E}$ can be identified with an algebraic $E_{8}\times E_{8}$-bundle on
$Z$ which is Mumford polystable with respect to the polarization
$\omega$. In view of this theorem one can use algebraic geometry to study the
moduli space of Heterotic vacua.

The type of bundles ${\mathcal E}$ allowed in a Heterotic
compactification is restricted in physics in three ways:
\begin{description}
\item[(Supersymmetry preservation)]  ${\mathcal E}$ has to be
Mumford polystable.
\item[(Anomaly cancellation)] $c_{2}({\mathcal E}) = c_{2}(Z)$.
\item[(Gauge symmetries)] If the compactification of the Heterotic
string  has
a group of symmetries $G \subset E_{8}\times E_{8}$, 
then the structure group of
${\mathcal E}$ can be reduced to the centralizer $G'$ of $G$ in
$E_{8}\times E_{8}$. Furthermore the corresponding $G'$ bundle
${\mathcal E}_{G'} \to Z$ should also be supersymmetric and anomaly-free.
\end{description}

\bigskip

Using these three principles one can look for special compactifications of the
Heterotic string that reproduce in their low energy limits well
understood and experimentally confirmed quantum field theories. Of
particular interest are compactifications that will lead to the
Standard Model of particle physics. For such compactifications one
imposes two additional requirements on the triple $(Z,{\mathcal E},\omega)$:
\begin{description}
\item[(Standard Model gauge symmetries)] The group $G$ of symmetries of
${\mathcal E}$, i.e. the centralizer inside $E_{8}\times E_{8}$
of a minimal subgroup $G' \subset E_{8}\times E_{8}$ to which the
structure group of ${\mathcal E}$ reduces, is
$G = U(1)\times SU(2)\times SU(3)$.
\item[($3$-generations condition)] 
\[
\chi(\op{ad}({\mathcal E})) = \frac{c_{3}(\op{ad}({\mathcal E}))}{2} = 3.
\]
\end{description}

Triples $(Z,\omega,{\mathcal E})$ satisfying the five conditions above
are hard to come by and tend to be rigid (see e.g. \cite{kachru}).
However the recent advances in string theory prompted by the
ground breaking work of Ho\v rava-Witten \cite{hw1,hw2,wittenSC} on orbifold
compactifications of M-theory, allow for a significant relaxation of
the anomaly cancellation condition. This leads to two essential
simplifications. First, it turns out that using the Ho\v rava-Witten
mechanisms one can suppress completely one copy of $E_{8}$ in the
structure group of ${\mathcal E}$. Secondly it was argued in
\cite{dolw,dow,ACK} that one can use M-theory 5-branes to relax the
equality in the anomaly cancellation condition to an inequality. This
leads to the following purely mathematical problem.

\

\bigskip

\noindent 
{\bf Main Problem.} Find a smooth Calabi-Yau 3-fold
$(Z,\omega)$ and a reductive subgroup $G' \subset E_{8}$ so that 
\begin{enumerate}
\item[$\lozenge$] the centralizer $G$ of $G'$ in $E_{8}$ is a group
isogenous to 
$SU(3)\times SU(2)\times U(1)$;
\item[$\lozenge$] there exists an $\omega$-stable $G'_{{\mathbb
 C}}$-bundle $\cV \to Z$ so that  
\begin{itemize}
\item $c_{1}(\cV) = 0$,
\item $c_{2}(Z) - c_{2}(\cV)$ is the class of an effective reduced curve
on $Z$,
\item $c_{3}(\cV) = 6$.
\end{itemize}
\end{enumerate}
\

\bigskip

\noindent
Here the Chern classes of $\cV$ are calculated in the adjoint
representation of $E_{8}$ considered as a representation of $G'$.
In fact for the physics applications it suffices for $G$ to contain a
group isogenous to $SU(3)\times SU(2)\times U(1)$ as a direct summand.

The groups $G' \subset E_{8}$ whose centralizer contains  $SU(3)\times
SU(2)\times U(1)$ as a direct summand 
can be classified. It turns out that there are no
connected subgroups $G'$ with $Z_{E_{8}}(G') = SU(3)\times
SU(2)\times U(1)$. The stability assumption on $\cV$ guarantees that the the
structure group of $\cV$ can not be reduced to a proper subgroup of
$G'_{{\mathbb C}}$. Therefore the structure group of the
associated $\pi_{0}(G'_{{\mathbb
C}})$-bundle $\cV\times_{G'_{\mathbb C}} \pi_{0}(G'_{\mathbb C})$ can
not be reduced to a proper subgroup of $\pi_{0}(G'_{\mathbb C})$. Since
$\cV\times_{G'_{\mathbb C}} \pi_{0}(G'_{\mathbb C})$ is a Galois
cover of $Z$ with Galois group $\pi_{0}(G'_{\mathbb C})$, this just
means that there should be a surjective homomorphism $\pi_{1}(Z)
\twoheadrightarrow \pi_{0}(G'_{\mathbb C})$ and so we are forced to
work with a  non-simply connected $Z$.

Some possible choices for $G'$ are: $SU(3)\times ({\mathbb Z}/6)$, 
$SU(4)\times ({\mathbb Z}/3)$ and $SU(5)\times ({\mathbb Z}/2)$. The
corresponding centralizers are isogenous to $(SU(3)\times
SU(2)\times U(1))\times U(1)\times U(1)$, $(SU(3)\times
SU(2)\times U(1))\times U(1)$ and $SU(3)\times
SU(2)\times U(1)$.
When $G'_{0}$, the connected component of the identity  in $G'$, is a classical
group,  it turns out that the Chern classes of $\cV$ in the 
fundamental representation of $G'_{0}$ coincide with the Chern classes
of $\cV$ in the adjoint representation of $E_{8}$.

In this paper we
explain how to build a big family of solutions of the {\bf Main
Problem} above for $G' = SU(5)\times ({\mathbb Z}/2)$. 

For concreteness we look for $Z$'s with $\pi_{1}(Z) = {\mathbb
Z}/2$. Let $\cV$ be an $SL(5,{\mathbb C})\times ({\mathbb Z}/2)$-bundle 
on such a $Z$. Then $\cV$ splits as a product of a rank five vector
bundle and the unique non-trivial local system on $Z$ with monodromy
${\mathbb Z}/2$. Pulling back this vector bundle to the universal
cover $X$ of $Z$ we get a rank five vector bundle on $X$ which is
invariant under the action of $\pi_{1}(Z)$ on $X$. Conversely every
$\pi_{1}(Z)$-equivariant vector bundle $V \to X$ descends to a vector
bundle on $Z$. Thus, in order to solve the {\bf Main Problem}, it
suffices to construct a quadruple $(X,\tau_{X},H,V)$ such that
the following conditions hold:
\

\bigskip

\begin{itemize}
\item[$({\mathbb Z}/2)$] $X$ is a smooth Calabi-Yau 3-fold and
$\tau_{X} : X \to X$ is a freely acting involution. $H$ is a fixed
ample line bundle (K\"{a}hler structure) on $X$.
\item[{\bf (S)}] $V$ is an $H$-stable vector bundle of rank five on $X$.
\item[{\bf (I)}] $V$ is $\tau_{X}$-invariant.
\item[{\bf (C1)}] $c_{1}(V) = 0$.
\item[{\bf (C2)}] $c_{2}(X) - c_{2}(V)$ is effective.
\item[{\bf (C3)}] $c_{3}(V) = 12$.
\end{itemize}
\

\bigskip

Since we need a mechanism for constructing bundles on $X$, we will
choose $X$ to be elliptically fibered and use the so called {\em
spectral construction} \cite{fmw,fmw-vb,donagi,bjps} to produce bundles
on $X$. Note that the spectral construction applies only to
elliptic fibrations, i.e. genus one fibrations with a section. This is
the reason we build an equivariant $V$ on $X$ rather than obtaining
directly $\cV$ on $Z$. In general, there are two ways in which the
spectral construction can be modified to work on genus one
fibrations such as $Z$. One is to work with a spectral cover in the
Jacobian fibration of $Z$ and an abelian gerbe on it. The other route
(which is the one we chose) is to work with equivariant spectral data
on the universal cover of $Z$. Note that there are higher algebraic
structures involved in both approaches: the stackiness of the first
approach is paralleled by  complicated group actions on the derived
category in the second.

More specifically we take $X$ to be a Calabi-Yau of Schoen
type \cite{schoenCY}, 
i.e. a fiber product of two rational elliptic surfaces $B$ and $B'$ over
${\mathbb P}^{1}$, both in the four dimensional family described in
\cite{dopw-i}. 
The rank five bundle $V$ is built as an extension of two
vector bundles $V_{2}$ and $V_{3}$ of ranks two and three
respectively. Each of these is manufactured by the spectral
construction. Alternatively $V$ may be viewed as a bundle
corresponding to spectral data with a {\em reducible} spectral
cover. Our preliminary research  of this problem (some
of which is recorded in \cite{usold}) showed that bundles
corresponding to smooth spectral covers are unlikely to satisfy all
of the above conditions. In fact, for the Calabi-Yau's we consider,
one can show rigorously (see Remark~\ref{rem-no-bundles}) 
that $V$'s coming from
smooth spectral covers can never satisfy {\bf (I)}, {\bf (C1)} 
 and {\bf (C3)} at the same time. 

The paper is organized as follows. 
Section~\ref{s-cy3} describes the construction of $X$ and lists
the geometric constraints on the spectral data which will ensure the
validity of {\bf (I)}. Section~\ref{s-spectral} deals with the actual
construction. We describe $V_{2}$ and $V_{3}$ in terms of their
spectral data. The data for each $V_{i}$ involves a spectral curve
$C_{i}$ in the surface $B$, a line bundle ${\mathcal N}_{i}$ on
$C_{i}$, another line bundle $L_{i}$ on the surface $B'$, and  
some optional parameters. The effect
of taking these additional parameters to be non-zero is 
interpretted in section~\ref{ss-hecke-patterns} as a series
of Hecke transforms. The freedom to perform these Hecke transforms
gives us at the end of the day infinitely many families of bundles. 
In section~\ref{ss-invariance} we explain how the geometric
information about the action of the spectral involution, 
obtained in \cite[Theorem~7.1]{dopw-i},   
takes care of condition {\bf (I)}. A delicate point here is that we
need two genericity assumptions on $C_{i}$. The first one is that
$C_{i}$ is finite over the base of the elliptic fibration on $B$. The
second assumption is that $\op{im}[\op{Pic}(B) \to \op{Pic}(C_{i})]$
is Zariski dense in $\op{Pic}^{0}(C_{i})$. In sections \ref{sss-C2}
and \ref{sss-C3} we check these two assumptions in the special case
that is untimately utilized in the construction of $V$.
In section~\ref{s-numerics}
we translate the remaining conditions into a sequence of rather
tight numerical inequalities. In Section~\ref{ss-solutions}  we show how
the latter can be
solved. In Section~\ref{s-summary} we summarize the construction of
$(X,\tau_{X},H,V)$ and give an estimate on the dimension of the moduli
space of 
such quadruples. Finally in Appendix~\ref{app-hecke} we have gathered
some basic facts on Hecke transforms of vector bundles which are used
in Section~\ref{s-spectral}.

\

\bigskip

{\bf Acknowledgements:} We would like to thank Ed Witten, Dima Orlov,
and Richard Thomas for valuable conversations on the subject of
this work. 

R.~Donagi is supported in part by an NSF grant DMS-9802456 as well as a
UPenn Research Foundation Grant. 
B.~A.~Ovrut is supported in part by a Senior Alexander von Humboldt
Award, by the DOE under contract No. DE-AC02-76-ER-03071 and by a
University of Pennsylvania Research Foundation Grant. 
T.~Pantev is supported in part by an NSF grant DMS-9800790 and by an
Alfred P. Sloan Research Fellowship. 
D.~Waldram would like to thank Enrico Fermi Institute at The
University of Chicago and the Physics Department of The Rockefeller
University for hospitality during the completion of this work.

\tableofcontents

\section{Elliptic Calabi-Yau threefolds with free involutions}
\label{s-cy3}

Our goal is to construct special $SU(5)$ bundles on
smooth Calabi-Yau 3-folds with fundamental group ${\mathbb Z}/2$. We
construct our Calabi-Yau 3-fold $Z$ as the quotient of a smooth
Calabi-Yau 3-fold $X$ by a freely acting involution $\tau_{X} : X \to
X$. Our $X$ will be elliptic and the elliptic fibration will be
preserved by $\tau_{X}$, so that $Z$ will still have a genus one
fibration. This enables us to apply the spectral
construction to produce bundles. 

The manifold $X$ is constructed as the fiber product $B\times_{\cp{1}}
B'$ of two rational elliptic surfaces $B$ and $B'$ which live in the
four dimensional family described in \cite[Section~4]{dopw-i}. 
For the first surface
$B$ we use the notation from \cite{dopw-i}. In
particular we have $\beta : B \to \cp{1}$, $e, \zeta : \cp{1} \to B$
and the involutions $\alpha_{B}, \tau_{B} : B \to B$ and
$\tau_{\cp{1}} : \cp{1} \to \cp{1}$. We use the same symbols
with primes for the corresponding objects on $B'$. Such constructions
were first considered in \cite{schoenCY}. In fact, exactly the
four dimensional subfamily of rational elliptic surfaces described in
\cite[Section~4]{dopw-i} happened to appear, in a different context, as an
example in \cite[Section~9]{schoenCY}.

We choose an isomorphism of $\cp{1}$ with ${\cp{1}}'$ which identifies
$\tau_{\cp{1}}$ with $\tau_{{\cp{1}}'}$ and sends $0 \in \cp{1}$ to
$\infty' \in {\cp{1}}'$ and $\infty \in \cp{1}$ to $0' \in
{\cp{1}}'$. With this convention we will make no
distinction between $\cp{1}$ and ${\cp{1}}'$ from now on.

Define $X := B\times_{\cp{1}} B'$. For a generic choice of $B$ and
$B'$ this $X$ will be smooth. It is
an elliptic 3-fold in two ways: via its projections $\pi : X \to B'$
and $\pi' : X \to B$. Since most of our analysis will involve the
elliptic fibers we will work with the elliptic structure $\pi : X \to
B'$ in order to avoid cumbersome notation. By construction the
discriminant of $\pi$ is in the linear system $\beta^{'*}{\mathcal
O}_{\cp{1}}(12) = -12K_{B'}$ and so $X$ is a Calabi-Yau 3-fold.

For the zero section of $\pi$ we take  the section 
$\sigma : B' \to X$ corresponding to $e :
\cp{1} \to B$.  Let $\alpha_{X} := \alpha_{B}\times_{\cp{1}} \tau_{B'}$
and let $\tau_{X} := \tau_{B}\times_{\cp{1}}
\tau_{B'}$. Since the fixed points of $\tau_{B}$ and $\tau_{B'}$ sit
over $\infty$ and $0$ respectively, we conclude that $\tau_{X}$ acts
freely on $X$. In particular the quotient $Z := X/\tau_{X}$ is
non-singular. We claim that $Z$ is in fact a Calabi-Yau. This is
equivalent to saying that $\tau_{X}$ preserves the holomorphic 3-form
on $X$. Indeed, $\tau_{X}$ acts on $H^{0}(X,\Omega^{3}_{X})$
as multiplication by a number $\lambda \in {\mathbb C}^{\times}$.  Since the
fiber $f_{0}\times f_{0}'$ of $X \to \cp{1}$ is stable under
$\tau_{X}$ it suffices to compute the action of $\tau_{X}$ on
$H^{0}(f_{0}\times f_{0}',\Omega^{3}_{X})$. But 
\[
\Omega^{3}_{X|f_{0}\times f_{0}'} = (K_{f_{0}}\boxtimes
K_{f_{0}'})\otimes T_{0}^{*}\cp{1},
\]
and $\tau_{X}$ acts as $\tau_{B|f_{0}}$ on $K_{f_{0}}$,
$\tau_{B'|f'_{0}}$ on $K_{f'_{0}}$ and as $\tau_{\cp{1}}$ on
$T_{0}^{*}\cp{1}$. Since  $\tau_{B|f_{0}}$ is a translation on
$f_{0}$, it acts on $H^{0}(f_{0}, K_{f_{0}})$ as $+1$. Since
$\tau_{B'|f'_{0}}$ and $\tau_{\cp{1}}$ each have a fixed point, they
act as $-1$ on $H^{0}(f'_{0},K_{f'_{0}})$ and $T_{0}^{*}\cp{1}$
respectively. Hence $\lambda = 1\cdot (-1)\cdot (-1) = 1$.

The vector bundles on $Z$ can be interpreted as
$\tau_{X}$-invariant vector bundles on $X$.  To construct vector
bundles on $X$ we will exploit the fact that $X$ is an elliptically
fibered 3-fold and so we can manufacture bundles by using a relative
Fourier-Mukai transform.

Concretely, let ${\mathcal P}_{X} \to X\times_{B'} X$ be the Poincare
sheaf corresponding to the section $\sigma$. That is ${\mathcal
P}_{X}$ is the rank one torsion-free sheaf given by
\[
{\mathcal P}_{X} = {\mathcal O}_{X\times_{B'} X}(\Delta -
\sigma\times_{B'} X - X\times_{B'}\sigma - m^{*}c_{1}(B')) =
p_{13}^{*}{\mathcal P}_{B},
\]
where $m : X\times_{B'} X \to B'$  and $p_{13} : X\times_{B'} X =
B\times_{\cp{1}} B' \times_{\cp{1}} B  \to B\times_{\cp{1}} B$ are 
the natural projections. As in \cite[Section~6]{dopw-i}, one argues that
${\mathcal P}_{X}$ defines an autoequivalence (see
\cite[Theorem~1.2]{bridgeland-cy}) of $D^{b}(X)$:
\[
\xymatrix@R=4pt{
\FM_{X} : & D^{b}(X) \ar[r] & D^{b}(X) \\
& {\mathcal F} \ar@{|->}[r] &
R^{\bullet}p_{1*}(p_{2}^{*}{\mathcal
F}\stackrel{L}{\otimes} {\mathcal P}_{X}).
}
\]
If $V \to X$ is a vector bundle of rank $r$ which is semistable and of
degree zero on each fiber of $\pi : X \to B'$, then its Fourier-Mukai
transform $\FM_{X}(V)[1]$ is a torsion sheaf of pure dimension two on
$X$. The support of $\FM_{X}(V)[1]$ is a 
surface $i_{\Sigma} : \Sigma \hookrightarrow X$ which is finite of 
degree $r$ over
$B'$. Furthermore $\FM_{X}(V)$ is of rank one on $\Sigma$. In fact, if
$\Sigma$  is
smooth, then $\FM_{X}(V)[1] = i_{\Sigma*}L$ is just the extension by
zero of some 
line bundle $L \in \op{Pic}(\Sigma)$. Conversely if ${\mathcal N} \to X$ is
a sheaf of pure dimension two which is flat over $B'$, then
$\FM_{X}({\mathcal N})$ is a vector bundle on $X$ of rank equal to the
degree of $\op{supp}({\mathcal N})$ over $B'$ and whose first Chern
class is vertical (for the projection $\pi : X \to B'$). This
correspondence between vector bundles on $X$ and sheaves on $X$
supported on finite covers of $B'$ is commonly known as the {\em spectral
construction} and has been extensively studied in the context of
Weierstrass elliptic fibrations \cite{fmw,fmw-vb,donagi,bjps}. The
torsion sheaf ${\mathcal N}$ on $X$ is called {\em spectral datum} and
the surface $\Sigma = \op{supp}{{\mathcal N}}$ is called a {\em spectral
cover}.

Since our elliptic Calabi-Yau $X$ is not Weierstrass we briefly
describe how the spectral construction works (at least for generic
spectral data) on $X$ and how it interacts with the involution
$\tau_{X}$. First we need to understand the action of $\FM_{X}$ on
line bundles on $X$. Note that since $X = B\times_{\cp{1}} B'$ is a
fiber product we have $\op{Pic}(X) = (\op{Pic}(B)\times
\op{Pic}(B'))/\op{Pic}(\cp{1})$. In particular, every line bundle on
$X$ can be written as $L\boxtimes L' := \pi^{'*}L\otimes \pi^{*}L'$
for some $L \to B$ and $L' \to B'$. 

\begin{lem} \label{lem-FMX-lb} For every line bundle ${\mathcal L} = 
L\boxtimes L'$ on $X$, the actions of the Fourier-Mukai transform and
of the spectral involution are given by:
\begin{itemize}
\item[{\em (a)}] $\FM_{X}({\mathcal L}) = 
\FM_{X}(L\boxtimes L')  = \pi^{'*}\FM_{B}(L)\otimes \pi^{*}L' =
\FM_{B}(L)\boxtimes L'$.
\item[{\em (b)}] $\T_{X}({\mathcal L}) := (\FM_{X}^{-1}\circ
\tau_{X}^{*} \circ
\FM_{X})({\mathcal L}) = \pi^{'*}(\T_{B}(L))\otimes
\pi^{*}(\tau_{B'}^{*}L') = \T_{B}(L)\boxtimes \tau_{B'}^{*}L'.$
\end{itemize}
\end{lem}
{\bf Proof.}  Part (b) is an obvious consequence of part (a). To prove
part (a) we will use the identification $X\times_{B'} X = B\times_{\cp{1}} B'
\times_{\cp{1}} B$. In terms of this identification we have:
\[
\begin{split}
\FM_{X}({\mathcal L}) & = Rp_{23*}(p_{12}^{*}{\mathcal L}\otimes
{\mathcal P}_{X}) \\
& = Rp_{23*}(p_{12}^{*}(L\boxtimes L')\otimes p_{13}^{*}{\mathcal
P}_{B}) \\
& = Rp_{23*}(\pi_{1}^{*}L\otimes \pi_{2}^{*}L'\otimes p_{13}^{*}{\mathcal
P}_{B}) \\
& = Rp_{23*}(p_{13}^{*}(p_{1}^{*}L\otimes {\mathcal
P}_{B})\otimes p_{23}^{*}(\pi^{*}L')) \\
& = Rp_{23*}(p_{13}^{*}(p_{1}^{*}L\otimes {\mathcal
P}_{B}))\otimes\pi^{*}L'.
\end{split}
\]
Here $\pi_{1}$, $\pi_{2}$ and $\pi_{3}$ are the natural projections of
$B\times_{\cp{1}} B' \times_{\cp{1}} B$ onto $B$, $B'$ and $B$
respectively, $p_{i} : B\times_{\cp{1}} B \to B$ are the
projections on the two factors, and in the last identity we have used
the projection formula for $p_{23}$.

Now using the base change property for the fiber square
\[
\xymatrix{
X\times_{B'} X \ar[r]^-{p_{23}} \ar[d]_-{p_{13}} & B'\times_{\cp{1}} B
\ar[d]^-{\pi'} \\
B\times_{\cp{1}} B \ar[r]_{p_{2}} & B
}
\]
we get $Rp_{23*}p_{13}^{*} = \pi^{'*}Rp_{2*}$ and so
\[
\FM_{X}({\mathcal L}) = \pi^{'*}(Rp_{2*}(p_{1}^{*}L\otimes {\mathcal
P}_{B}))\otimes \pi^{*}L' = \FM_{B}(L)\boxtimes L'.
\]
The lemma is proven. \hfill $\Box$

Let now $i_{\Sigma} : \Sigma \hookrightarrow X$ be a surface which is
finite and of 
degree $r$ over $B'$. Then for any line bundle ${\mathcal L} \in \op{Pic}(X)$
the torsion sheaf ${\mathcal N} := i_{\Sigma*}i^{*}_{\Sigma}{\mathcal
L}$ has a
resolution  by global line bundles. Namely 
\[
0 \to {\mathcal L}(-\Sigma) \to {\mathcal L} \to {\mathcal N} \to 0.
\]
In particular ${\mathcal N}$ is quasi-isomorphic to the two-step
complex of line bundles $[{\mathcal L}(-\Sigma) \to {\mathcal L}]$ on $X$
and so the actions of $\FM_{X}$ and $\T_{X}$ on ${\mathcal N}$ can be
computed via the formulas in Lemma~\ref{lem-FMX-lb}. Specifically we
have:

\begin{lem} \label{lem-generic-tauX-invariance} Let ${\mathcal L} =
L\boxtimes L'$ be a global line bundle on $X$ and let $i_{\Sigma} : \Sigma
\hookrightarrow X$ be a surface finite over $B'$. Let ${\mathcal N} =
i_{\Sigma*}i^{*}_{\Sigma}{\mathcal L}$ be such that 
$V = \FM_{X}({\mathcal N})$ is a
rank $r$ vector bundle on $X$ with $c_{1}(V) = 0$. Then $\tau_{X}^{*}V
\cong V$ if and only if the following three conditions
\begin{itemize}
\item $\alpha_{X}(\Sigma) = \Sigma$;
\item $\tau_{B'}^{*}L' \cong L'$;
\item $\T_{B}^{*}L \cong L$.
\end{itemize}
are satisfied.
\end{lem}

\bigskip

\

\begin{rem} \label{rem-no-bundles}
Notice that the $\tau_{X}$ invariance of $V$ amounts to two separate
conditions on the spectral data. The first is that the spectral
surface $\Sigma$ has to be invariant under the involution
$\alpha_{X}$. This condition is relatively
easy to satisfy. It just means that $\Sigma$ is pulled back from the
quotient $X/\alpha_{X}$. The second condition requires the $\tau_{B'}$
invariance of $L'$ and the $\T_{B}$ invariance of $L$. 

In fact, the formulas in \cite[Table~3]{dopw-i} (written in terms of
the basis of $H^{2}(B,{\mathbb Z})$ described in
\cite[Section~4.2]{dopw-i})  show that $L \in
\op{Pic}(B)\otimes {\mathbb Q}$ will be $\T_{B}$-invariant if and only if
$L$ is in the affine subspace 
\[
- \frac{1}{2} e_{1} + \op{Span}_{{\mathbb Q}}(f, e_{9}, e_{4} - e_{5},
  e_{4} - e_{6}, 3\ell - 2(e_{4} + e_{5} + e_{6}) - 3e_{7}, \ell -
  e_{7} - 2e_{8}),
\]
which does not intersect $\op{Pic}(B) \subset \op{Pic}(B)\otimes
{\mathbb Q}$. This implies that $V$ can not be $\tau_{X}$-invariant if
${\mathcal N} = i_{\Sigma*}i^{*}_{\Sigma}{\mathcal L}$ for some global
${\mathcal L} \in \op{Pic}(X)$. For $\Sigma$ smooth and very ample the
Lefschetz hyperplane section theorem asserts that every ${\mathcal N}$
comes from a global ${\mathcal L}$ and hence one is forced to work
with singular or non-very ample surfaces $\Sigma$.
\end{rem}

\section{The construction} \label{s-spectral}

\subsection{The basic construction} \label{ss-basic}

In this section we describe in detail our method of constructing
$\tau_{X}$-invariant vector bundles on $X$.

In order to circumvent the difficulty pointed out in Remark
\ref{rem-no-bundles} we build our rank five bundle $V$ on $X$ not
directly by the spectral construction but as an extension
\[
0 \to V_{2} \to V \to V_{3} \to 0.
\]
Here $V_{i}$, $i =2,3$ is a rank $i$ bundle on $X$ which is
$\tau_{X}$-invariant and satisfies some strong numerical conditions
which will be discussed in the next section. In addition, we will see
that the stability condition on $V$ amounts to the extension being
non split.

Each $V_{i}$ is produced by an application of the spectral construction
on $X$ with a reducible spectral cover and a
a line bundle on it which is {\em not} the restriction of a global line bundle
on $X$. Define $V_{i}$ from its spectral data as follows:
\begin{itemize}
\item Let $C_{i}$ be a  curve in the linear system $|{\mathcal
O}_{B}(ie+k_{i}f)|$ where $k_{i}$ is an integer.  Let $\Sigma_{i} :=
C_{i}\times_{\cp{1}} B'$. Recall that $\beta' : B' \to \cp{1}$ has two
$I_{2}$ fibers $f'_{1}, f'_{2}$. Let $F_j, j=1,2$ be the corresponding
fibers of $B$; note that, because of the way we glued the $\cp{1}$
bases, these are {\em not} the reducible fibers $f_j$ of $B$.
Let 
\[
\{ p_{ijk} \}_{k = 1}^{i} := C_{i} \cap F_{j}.
\]
Then $\Sigma_{i} \to C_{i}$ is an elliptic
surface having $2i$ fibers of type $I_{2}$: $(n_{j}'\cup o_{j}')\times
\{ p_{ijk}\}$ where $j = 1, 2$ and  $k = 1, \ldots, i$.
\item Define $V_{i}$ as
\[
V_{i} = \FM_{X}\left(\left(\Sigma_{i}, (\pi'_{|\Sigma_{i}})^{*}
{\mathcal N}_{i}\otimes
\pi^{*}L_{i}\otimes {\mathcal O}_{\Sigma_{i}}\left(-\sum \{ p_{ijk}
\}\times ( a_{ijk}
n_{j}' + b_{ijk} o_{j}')\right)\right)\right),
\]
where $L_{i}$ is a line bundle on $B'$, ${\mathcal N}_{i}$ is a line
bundle on the curve $C_{i}$ and the optional parameters 
$a_{ijk},b_{ijk}$ are integers.
\end{itemize}

\bigskip 
Note that there is a redundancy in our choices, because $n'_j+o'_j=f'_j$
is a pullback from $\cp{1}$, and so can be absorbed in $L_i$. In
particular, we can always arrange for all the coefficients
$a_{ijk},b_{ijk}$ to be non-negative. Also, without a loss of
generality we may assume that for any given $j$ we have $a_{ijk}\cdot
b_{ijk} = 0$ for all $i,k$. With this convention, we have an
alternative description of $V_i$: Put
\[
\widetilde{W}_{i} := V_{i}\otimes \pi^{*} L_{i}^{-1}. 
\]
It turns out that the bundle 
$\widetilde{W}_{i}$ can be constructed directly. Consider the vector
bundle $W_{i}$ on $B$, built by the spectral construction
as $W_{i} := \FM_{B}(C_{i},{\mathcal N}_{i})$. Then
$\widetilde{W_{i}}$ is obtained from the vector bundle
$\pi^{'*}W_{i}$ by $a_{ijk}$ successive Hecke transforms along the
divisors $\{ p_{ijk} \}\times n_{j}'$ and $b_{ijk}$ successive Hecke
transforms  along the
divisors $\{ p_{ijk} \} \times o_{j}'$.  The
center of each Hecke is a 
line bundle on the surface $F_{j}\times n_{j}'$ or $F_{j}\times
o_{j}'$. (For the definition and basic properties of Hecke transforms 
see appendix~\ref{app-hecke}).

In fact, $V$ itself could be built by applying the
spectral construction on $X$ to the reducible spectral cover
$\Sigma_{2} \cup \Sigma_{3}$ and an appropriately chosen sheaf on it.
However the construction with extensions is technically
easier because it allows us to avoid dealing with sheaves on singular
surfaces. This approach is a variation on the method employed by Richard
Thomas in \cite{richardt}.

\begin{rem} \label{rem-freedom in-spectral-cover} Observe that in the
definition of $C_{i}$ we could have taken the linear system more
generally to be of the form $|{\mathcal
O}_{B}(ie+k_{i}f + \eta_{i})|$ where $k_{i}$ is an integer and
$\eta_{i} \in \op{Pic}(B)$ is a class perpendicular to $e$ and $f$. If we
impose the condition $c_{1}(W_{i}) = 0$, then the classes $\eta_{i}$
are forced to be zero by the Riemann-Roch formula. However the introduction of
the $L_{i}$'s gives us the extra freedom of working with $W_{i}$'s
that have arbitrary vertical $c_{1}$. We will not exploit this extra
freedom but we expect that many examples exist which are similar to
ours but have $\eta_{i} \neq 0$.
\end{rem}

\

\medskip

\noindent
Since the Hecke interpretation of $\widetilde{W}_{i}$ 
will be important in determining the
invariance properties of $V$ and in implementing the numerical
constraints, we proceed to spell it out explicitly.

\subsection{Reinterpretation via Hecke transforms} \label{ss-hecke-patterns}

Recall from section~\ref{s-cy3} that $X = B\times_{\cp{1}} B'$ fits
into a commutative diagram of projections
\begin{equation} \label{eq-main-diagram}
\xymatrix{
& X \ar[dl]_-{\pi'} \ar[dr]^-{\pi} & \\
B \ar[dr]_-{\beta} & & B' \ar[dl]^-{\beta'} \\
& \cp{1} &
}
\end{equation}
Let $C \subset B$ be a smooth connected curve in the linear system
$|{\mathcal O}_{B}(re + kf)|$. Let ${\mathcal N} \in \op{Pic}^{d}(C)$
and let $W := \FM_{B}((C,{\mathcal N}))$ be the corresponding rank $r$
vector bundle on $B$.

Consider $\Sigma = \pi^{'-1}(C) = C\times_{\cp{1}} B'$ and the line bundle
${\mathcal L} = (\pi'_{|\Sigma})^{*}{\mathcal N} = {\mathcal
N}\boxtimes {\mathcal O}_{B'}$ on $\Sigma$.

Let $f_{j}' = n_{j}'\cup o_{j}'$, $F_{j} =
\beta^{-1}(\beta'(f_{j}'))$, $j = 1, 2$
be as above. We assume that $C$ is general enough so that the
intersections $C\cap F_{j}$, $j = 1,2$ are transversal. 

Let $p \in C\cap F_{j}$ and let $a$ be a non-negative integer. Define
\[
W[a,p] := \FM_{X}((\Sigma, {\mathcal L}(- a(\{ p \}\times
n_{j}')))).
\]
Consider the divisor $D = F_{j}\times n_{j}' \subset X$ and the line
bundles
\[
{\mathcal O}_{F_{j}}(p-e)\boxtimes {\mathcal
O}_{n_{j}'}(2a) \in \op{Pic}(D),
\]
where $a \in {\mathbb Z}$ and by a abuse of notation $e$ denotes the
point of intersection  of the curves $e, F_{j}$ in
$B$. With this notation we have
\begin{lem} \label{lem-hecke-picture} Fix $a \geq 0$.  
\begin{itemize}
\item[{\bf (i)}] There is a canonical surjective map
$W[a,p]_{|D} \to 
{\mathcal O}_{F_{j}}(p-e)\boxtimes {\mathcal 
O}_{n_{j}'}(2a)$ which fits in a short exact sequence of vector
bundles on $D$
\begin{equation*}
0 \to K_{a} \to W[a,p]_{|D} \to
{\mathcal O}_{F_{j}}(p-e)\boxtimes {\mathcal
O}_{n_{j}'}(2a) \to 0. \tag{$\psi_{a+1}$}
\end{equation*}
\item[{\bf (ii)}] For $a=0$ we have $W[0,p] = \pi^{'*}W$ and
for $a \geq 1$
\[
W[a,p] = \heck^{-}_{(\psi_{a})}\circ
\heck^{-}_{(\psi_{a-1})}\circ \ldots \circ \heck^{-}_{(\psi_{1})}(\pi^{'*}W).
\]
\end{itemize}
\end{lem}
{\bf Proof.} We will prove the lemma by induction in $a$. By
definition we have $W[0,p] = \pi^{'*}W$ which takes care of
the base of the induction. Assume that $W[i,p] =
\hdown{\psi_{i}}{W[i-1,p]}$ for all $0< i < a$. Consider the
short exact sequence of sheaves on $\Sigma$:
\begin{equation} \label{eq-L-ses}
0 \to {\mathcal L}(-a(\{p\}\times n_{j}')) \to {\mathcal
L}(-(a-1)(\{p\}\times n_{j}')) \to {\mathcal
L}(-(a-1)(\{p\}\times n_{j}'))_{|\{p\}\times n_{j}'} \to 0.
\end{equation}
We have ${\mathcal L}_{|\{p \}\times n_{j}'} =
((\pi'_{|\Sigma})^{*}{\mathcal N})_{|\{p \}\times n_{j}'} = {\mathcal
O}_{\{p \}\times n_{j}'}$. Also $\{ p \}\times n_{j}'$ is a component
of an $I_{2}$ fiber of the elliptic surface $\pi'_{|\Sigma} : \Sigma
\to C$ and so ${\mathcal O}_{\{ p \}\times n_{j}'}(\{ p \}\times
n_{j}') = {\mathcal O}_{\{ p \}\times n_{j}'}(-2)$. Let $i_{\Sigma} :
\Sigma \hookrightarrow X$ and $\iota : n_{j}' = \{ p \}\times n_{j}'
\hookrightarrow \Sigma \subset X$ denote the natural inclusions. Then
if we extend each of the sheaves in the sequence \eqref{eq-L-ses} by
zero we obtain a short exact sequence of sheaves on $X$:
\begin{equation} \label{eq-L-FM-ses}
0 \to i_{\Sigma*}{\mathcal L}(-a(\{p\}\times n_{j}')) \to i_{\Sigma*}{\mathcal
L}(-(a-1)(\{p\}\times n_{j}')) \to \iota_{*}{\mathcal
O}_{n_{j}'}(2(a-1)) \to 0.
\end{equation}
Applying $\FM_{X}$ to \eqref{eq-L-FM-ses} we get 
\begin{equation} \label{eq-Wa-hecke}
0 \to W[a,p] \to W[a-1,p] \to
\FM_{X}(\iota_{*}{\mathcal 
O}_{n_{j}'}(2(a-1))) \to 0.
\end{equation}
By the definition of $\FM_{X}$ we have 
\[
\FM_{X}(\iota_{*}{\mathcal O}_{n_{j}'}(2(a-1))) =
Rp_{2*}^{\bullet}(p_{1}^{*}(\iota_{*}{\mathcal
O}_{n_{j}'}(2(a-1))\stackrel{L}{\otimes} {\mathcal P}_{X}),
\]
where $p_{1}, p_{2} : X\times_{B'} X \to X$ are the natural
projections. 

If we use the identification $X\times_{B'} X \cong
B\times_{\cp{1}}B\times_{\cp{1}} B'$, then the projection $p_{i} :
X\times_{B'} X \to X$ gets identified with the projection $p_{i3} :
B\times_{\cp{1}}B\times_{\cp{1}} B' \to B\times_{\cp{1}} B'$ and
${\mathcal P}_{X} = p_{12}^{*}{\mathcal P}_{B}$. In particular in
terms of
the identification $X\times_{B'} X \cong
B\times_{\cp{1}}B\times_{\cp{1}} B'$ we see that
$p_{1}^{*}(\iota_{*}{\mathcal O}_{n_{j}'}(2(a-1)))$ is supported on
the surface $\{ p \}\times F_{j}\times n_{j}' \subset
B\times_{\cp{1}}B\times_{\cp{1}} B'$ is precisely
\[
\begin{split}
\op{pr}_{n_{j}'}^{*}{\mathcal O}_{n_{j}'}(2(a-1))\otimes {\mathcal
P}_{X|\{ p \}\times F_{j} \times n_{j}'} & =
\op{pr}_{n_{j}'}^{*}{\mathcal O}_{n_{j}'}(2(a-1))\otimes
\op{pr}_{F_{j}}^{*}{\mathcal P}_{B|\{p \}\times F_{j}} \\
& =
\op{pr}_{n_{j}'}^{*}{\mathcal O}_{n_{j}'}(2(a-1))\otimes
\op{pr}_{F_{j}}^{*}{\mathcal O}_{n_{j}'}(p - e).
\end{split}
\] 
Also for the restricted map $p_{2|\{ p \}\times F_{j} \times n_{j}'} :
\{ p \}\times F_{j} \times n_{j}' \to X$ we get 
\[
p_{2|\{ p \}\times F_{j} \times n_{j}'} = p_{23|\{ p \}\times F_{j}
\times n_{j}'} = i_{D},
\]
and hence
\[
\begin{split}
\FM_{X}(\iota_{*}{\mathcal O}_{n_{j}'}(2(a-1))) & =
R_{p_{2*}}^{\bullet}(p_{1}^{*}(\iota_{*}{\mathcal
O}_{n_{j}'}(2(a-1))\stackrel{L}{\otimes} {\mathcal P}_{X}) \\
& =
i_{D*}({\mathcal O}_{F_{j}}(p-e)\boxtimes {\mathcal O}_{n_{j}'}(2(a-1))).
\end{split} 
\]
which combined with \eqref{eq-Wa-hecke} 
concludes the proof of the lemma. \hfill $\Box$

\bigskip

\noindent
If now $b \geq 0$ is another integer we may consider also the vector bundle
\[
W\{b,p\} = \FM_{X}((\Sigma, {\mathcal L}(-b(\{p\}\times o_{j}')))).
\]
In exactly the same way we see that $W\{0\} = \pi^{'*}W$,
that for every $b \geq 1$ there is a canonical exact sequence
\begin{equation*}
0 \to M_{b} \to W\{b,p\}_{|\{F_{j} \}\times o_{j}'} \to {\mathcal
O}_{F_{j}}(p-e)\boxtimes {\mathcal O}_{o_{j}'}(2b) \to 0, \tag{$\phi_{b+1}$} 
\end{equation*} 
and that 
\[
W\{b,p\} = \heck^{-}_{(\phi_{a})}\circ
\heck^{-}_{(\phi_{a-1})}\circ \ldots \circ \heck^{-}_{(\phi_{1})}(\pi^{'*}W).
\]
For future reference we record
\begin{cor} \label{cor-chern} The Chern classes of $W[a,p]$ and
$W\{b,p\}$ are given by
\[
\begin{split}
ch(W[a,p]) & = \pi^{'*}ch(W) - a\pi^{*}n_{j}' - a^{2}(f\times \op{pt});
\\
ch(W\{b,p\}) & = \pi^{'*}ch(W) - b\pi^{*}o_{j}' - b^{2}(f\times \op{pt})
\end{split}
\]
\end{cor}
{\bf Proof.} Clearly it suffices to prove the corollary for
$W[a,p]$. By Lemma~\ref{lem-hecke-picture} we have short exact
sequences
\[
0 \to W[n,p] \to W[n-1,p] \to i_{D*}\psi_{n} \to 0
\]
for all $n \geq 1$. Here we have slightly abused the notation
by writing $\psi_{a}$ for the middle term of the short exact sequence
$(\psi_{a})$. Hence $ch(W[n,p]) = ch(W[n-1,p]) -
ch(i_{D*}\psi_{n})$ and so
\[
ch(W[a,p]) = \pi^{'*}ch(W) - \sum_{n=1}^{a} ch(i_{D*}\psi_{n}).
\]
Using Grothendieck-Riemann-Roch we calculate
\[
\begin{split}
ch(i_{D*}\psi_{n}) & = i_{D*}(ch(\psi_{n})td(D))td(X)^{-1} \\
& = i_{D*}\left(\left(1 + \psi_{n} +
\frac{\psi_{n}^{2}}{2}\right)(1 + F_{j}\times \op{pt})\right)(1 -
(f\times \op{pt} + \op{pt}\times f')) \\
& = i_{D*}((1 + 2(n-1)F_{j}\times \op{pt})(1 + F_{j}\times
\op{pt}))(1 - (f\times \op{pt} + \op{pt}\times f')) \\
& = i_{D*}(1 + (2n-1)F_{j}\times \op{pt})(1 - (f\times \op{pt} +
\op{pt}\times f')) \\
& = D + (2n-1)(f\times \op{pt}) = \pi^{*}n_{j}' + (2n-1)(f\times \op{pt}).
\end{split} 
\]
Consequently
\[
\begin{split}
ch(W[a,p]) & = \pi^{'*}ch(W) - \sum_{n = 1}^{a}(\pi^{*}n_{j}' +
(2n-1)(f\times \op{pt}) \\
& = \pi^{'*}ch(W) - a\pi^{*}n_{j}' - a^{2}(f\times \op{pt}).
\end{split}
\]
The corollary is proven.  \hfill $\Box$

\bigskip

\noindent
Finally, we are ready to give the Hecke interpretation of 
$\widetilde{W}_{i} = V_{i}\otimes \pi^{*}L_{i}^{-1}$. Recall that 
\[
\widetilde{W}_{i} = \FM_{X}\left(\left(\Sigma_{i},{\mathcal
L}_{i}\left(-\sum 
\{ p_{ijk}\}\times (a_{ijk}n_{j}' + b_{ijk}o_{j}')\right)\right)\right)
\]
where $a_{ijk}, b_{ijk}$ are non-negative integers satisfying
$a_{ijk}b_{ijk} =0$. Since Hecke transforms whose centers have
disjoint supports obviously commute, we see from the above discussion
that 
\begin{equation} \label{eq-hecke-pattern}
\widetilde{W}_{i} =
W[a_{i11},p_{i11}]\{b_{i11},p_{i11}\}[a_{i12},p_{i12}]\{b_{i12},p_{i12}\}
\ldots [a_{i1i},p_{i1i}]\{b_{i1i},p_{i1i}\}.
\end{equation}

\section{Invariant spectral data} \label{ss-invariance}

In this section we examine the conditions for $V$ to be
$\tau_{X}$-invariant.  It is easy to reduce this, first to invariance
of the $V_{i}$, then to invariance of the $W_{i}$. Indeed, assume
that the bundles $V_i$ are
$\tau_{X}$-invariant, and choose liftings of the $\tau_{X}$ action
to the $V_i$. The space $\op{Ext}^1(V_3,V_2)$
parameterizing all extensions is a direct sum of its invariant and
anti-invariant subspaces. So if $\op{Ext}^1(V_3,V_2) \neq 0$ we also have 
an extension which is either invariant or anti-invariant. Finally,
changing the lifted action of $\tau_X$ on one of the $V_i$ interchanges
invariants with anti-invariants, so we are done.

Since $V_{i} = \widetilde{W}_{i}\otimes \pi^{*} L_{i}$, we have
\[
\tau_X^*(V_i) = \tau_X^*(\widetilde{W}_{i}) \otimes  \pi^{*}\tau_{B'}^* L_{i}.
\]
So it suffices to have  a
$\widetilde{W}_{i}$ which is $\tau_X$-invariant and an $L_i$ which is
$\tau_{B'}$-invariant. From \cite[Table~1]{dopw-i} we know that there is a
$6$-dimensional lattice of $\tau_{B'}$-invariant classes on $B'$, so we
have lots of possibilities for the $L_i$. Now
$\widetilde{W}_{i}$ is a Hecke transform of $\pi^{'*}(W_i)$, so we
want $W_i$ to be $\tau_B$-invariant and the center of the Hecke
transform to be $\tau_X$-invariant. Above we took the support of this
Hecke to be an arbitrary collection of components of the surfaces $F_j
\times f_j'$ for $j=1,2$. It can
be seen from \cite[Table~1]{dopw-i} together with
the expression \cite[Formula~(4.2)]{dopw-i} for the components of the
$I_{2}$-fibers of $B$ in terms of our basis,  that the action of $\tau_{B'}$
interchanges $o_{1}$ with $n_{2}$ and $o_{2}$ with $n_{1}$. Therefore
the condition for $\tau_{X}$-invariance of the center of the Hecke
transform becomes 
$a_{i1k} = b_{i2k}$ and $a_{i2k} = b_{i1k}$. Because of the redundancy
in our choices we are free to take $a_{i2k} = b_{i1k} = 0$ and
$a_{i1k} = b_{i2k} \geq 0$. 

Finally we have to find the conditions that will ensure the
$\tau_{B}$-invariance of $W_{i}$. 

\

\bigskip

\subsection{The $\tau_{B}$-invariance of $W_{i}$} 

Throughout
this subsection we work with a spectral curve $C_{i}$ in the linear system
$|ie + k_{i}f|$, $i = 2$ or $3$, which is {\em finite} over $\cp{1}$,
and a line bundle ${\mathcal N}_{i} \in \op{Pic}(C_{i})$.  
The $\tau_{B}$-invariance of $W_{i} =
\FM_{B}((C_{i},{\mathcal N}_{i}))$ is equivalent to the $T_{B}$
invariance of  $i_{C_{i}*}{\mathcal N}_{i}$. 

In 
\cite[Proposition~7.7]{dopw-i} we saw  that for any curve $C
\subset B$ which is finite over $\cp{1}$ and for any line bundle
${\mathcal N}$ on $B$ the image 
$\T_{B}(i_{C*}{\mathcal N})$ is again a sheaf of the form
$i_{\alpha_{B}(C)*}(?)$ 
for some line bundle $? \in \op{Pic}(\alpha_{B}(C))$. 
Therefore $\T_{B}$ induces a
well defined map $\T_{C} : \op{Pic}(C) \to
\op{Pic}(\alpha_{B}(C))$. Due to this,  
the $\tau_{B}$-invariance of $\FM_{B}((C,{\mathcal N}))$ is equivalent
to the following 
two conditions:
\begin{align}
C     & = \alpha_{B}(C) \label{eq-ci-inv}\\
{\mathcal N}     & = \T_{C    }({\mathcal N}). \label{eq-ni-invariance}
\end{align}

\begin{lem} \label{lem-vertical} The linear system $|re + kf|$
contains smooth $\alpha_{B}$-invariant curves if $r = 3$, $k \geq 3$
or if $r = 2$ and $k \geq 2$ is even.
\end{lem}
{\bf Proof.} First of all, from the explicit equations of a spectral
curve \cite{fmw} and Bertini's theorem, it is easy to see that the
general curve $C$  in the linear system $|re + kf|$ will 
be smooth as long as $k \geq r > 1$.  The same kind of analysis allows
one to understand the $\alpha_{B}$-invariant members of these linear
systems as well. Indeed, recall (see
e.g. \cite{fmw-vb}) that 
for every $a \geq 0$ we have an isomorphism 
$\beta_{*}{\mathcal O}_{B}(a e) = {\mathcal O}_{\cp{1}} \oplus
{\mathcal O}_{\cp{1}}(-2)\oplus \ldots {\mathcal O}_{\cp{1}}(-a)$. In
particular, by the projection formula we get isomorphisms
\[
\begin{split}
H^{0}(B,{\mathcal O}_{B}(e)) & = H^{0}(\cp{1},{\mathcal O}_{\cp{1}}) \\
H^{0}(B,{\mathcal O}_{B}(2e + 2f)) & = H^{0}(\cp{1},{\mathcal
O}_{\cp{1}}(2))\oplus H^{0}(\cp{1},{\mathcal
O}_{\cp{1}}) \\
H^{0}(B,{\mathcal O}_{B}(3e + 3f)) & = H^{0}(\cp{1},{\mathcal
O}_{\cp{1}}(3))\oplus H^{0}(\cp{1},{\mathcal
O}_{\cp{1}}(1)) \oplus H^{0}(\cp{1},{\mathcal
O}_{\cp{1}}).
\end{split}
\]
Let $X \in H^{0}(B,{\mathcal O}_{B}(2e + 2f))$, $Y \in
H^{0}(B,{\mathcal O}_{B}(3e + 3f))$ and $Z \in H^{0}(B,{\mathcal
O}_{B}(e)$ be the preferred sections corresponding to the generator of
the piece $ H^{0}(\cp{1},{\mathcal O}_{\cp{1}})$ under the above
decompositions. Note that in terms of the sections $x \in
H^{0}(P,{\mathcal O}_{P}(1)\otimes p^{*}{\mathcal O}_{\cp{1}}(2)$, $y \in
H^{0}(P,{\mathcal O}_{P}(1)\otimes p^{*}{\mathcal O}_{\cp{1}}(3)$ and
$z \in H^{0}(P,{\mathcal O}_{P}(1)$, which were used in
\cite[Section~3.2]{dopw-i} to define the Weierstrass model of $B$,
we have $x_{|W_{\beta}} = XZ$, $y_{|W_{\beta}} = Y$, $z_{|W_{\beta}} = Z^{3}$.

With this notation the isomorphism 
\[
H^{0}({\mathcal
O}_{\cp{1}}(k))\oplus H^{0}({\mathcal O}_{\cp{1}}(k-2)\oplus \ldots
\oplus H^{0}({\mathcal O}_{\cp{1}}(k-r) \to H^{0}(B,{\mathcal
O}_{B}(re + kf))
\]is given explicitly by the formula 
\[
(a_{k},a_{k-2},
\ldots, a_{k-r}) \mapsto (\beta^{*}a_{k})Z^{r} + (\beta^{*}a_{k-2})XZ^{r-2} +
(\beta^{*}a_{k-3})YZ^{r-3} + \ldots .
\] 
In particular the curves $C_{2}$ and $C_{3}$ can be identified with the
zero loci of 
\begin{equation} \label{eq-C2andC3}
(\beta^{*}a_{k_{2}})Z^{2} + (\beta^{*}a_{k_{2}-2})X \quad  \text{and}  \quad
(\beta^{*}a_{k_{3}})Z^{3} + (\beta^{*}a_{k_{3}-2})XZ +
(\beta^{*}a_{k_{3}-3})Y,
\end{equation}
respectively. 

Now recall, that in \cite[Section~3.2]{dopw-i} we identified
$\alpha_{B}$ with the involution induced from $\tau_{P|W_{\beta}}$ and
that $\tau_{P}$ acts trivially on the sections $x$, $y$, $z$. In view
of the comparison formulas  $x_{|W_{\beta}} = XZ$, $y_{|W_{\beta}} = Y$,
$z_{|W_{\beta}} = Z^{3}$, this implies that $\alpha_{B}^{*}(X) = X$, 
$\alpha_{B}^{*}(Y) = Y$ and $\alpha_{B}^{*}(Z) = Z$. Here the lifting
of the action of $\alpha_{B}$ to an action on line bundles of the form
${\mathcal O}_{B}(re + kf)$ is chosen in the way described in
\cite[Section~3.2]{dopw-i}. In particular
$\alpha_{B}^{*}(\beta^{*}s) = \beta^{*}(\tau_{\cp{1}}^{*}s)$ for any
section $s \in {\mathcal O}_{\cp{1}}(k)$. 

Since $\alpha_{B}^{*}$ acts linearly on the projective space $|re +
kf|$ it follows  that $\alpha_{B}$ will preserve a divisor $C \in |re +
kf|$ if and only if 
\[
C \in {\mathbb P}(H^{0}(B,{\mathcal O}_{B}(re +
kf))^{+})\cup {\mathbb P}(H^{0}(B,{\mathcal O}_{B}(re +
kf))^{-}) \subset |re + kf|,
\] 
where $H^{0}(B,{\mathcal O}_{B}(re +
kf))^{\pm}$ denote the $\pm 1$ eigenspaces of $\alpha_{B}^{*}$ acting
on $H^{0}(B,{\mathcal O}_{B}(re + kf))$. Therefore we see that for
each $i$ there are two families of $\alpha_{B}$-invariant $C_{i}$'s,
each parameterized by a projective space. In particular we will have
$\alpha_{B}(C_{i}) = C_{i}$ if and only if all the coefficients in the
polynomial expressions \eqref{eq-C2andC3} are simultaneously
$\tau_{\cp{1}}$-invariant or simulatneously $\tau_{\cp{1}}$-anti-invariant.
Now the Bertini theorem immediately implies that we can find a smooth
$C_{3}$, which is preserved by $\alpha_{B}$ as long as $k_{3} \geq 3$
and we can find a smooth
$C_{2}$, which is preserved by $\alpha_{B}$ as long as $k_{2} \geq 2$
and $k_{2}$ is even.  \hfill $\Box$ 

\bigskip

\

\begin{rem} \label{rem-vertical}
Unfortunately, when $k_{2}$ is odd the linear systems $|2e +
k_{2}f|^{\pm}$ will each have a fixed component and so all the
$\alpha_{B}$-invariant curves $C_{2}$ will be reducible. In order to
see this consider the homogeneous coordinates
$(t_{0}:t_{1})$ on $\cp{1}$ which were used
in \cite[Section~3.2]{dopw-i} to define the standard
action of $\tau_{\cp{1}}$. In other words $(t_{0}:t_{1})$ are such
that $\tau_{\cp{1}}^{*}(t_{0}) = t_{0}$, $\tau_{\cp{1}}^{*}t_{1} = -
t_{1}$ and $0 = (1:0)$ and $\infty = (0:1)$. Now it is clear that if
$a$ is a $\tau_{\cp{1}}$-invariant homogeneous polynomial in $t_{0}$
and $t_{1}$ of odd degree, then $a$ is divisible by $t_{0}$. Similarly
if $a$ is a $\tau_{\cp{1}}$-anti-invariant polynomial of odd degree,
then $a$ is divisible by $t_{1}$. In particular, since $k_{2}$ and
$k_{2} - 2$ have the same parity we see that for $k_{2}$ odd, the
fiber $f_{\infty}$ is the  fixed component of the
linear system  $|2e + k_{2}f|^{+}$ and the fiber $f_{0}$ is the fixed
component of the linear system  $|2e + k_{2}f|^{-}$.
\end{rem}

\

\begin{lem} \label{lem-invariant-mod-density} Let $C$ be an
$\alpha_{B}$-invariant curve in $|re + kf|$ which is finite over
$\cp{1}$, and assume that $i_{C}^{*}\op{Pic}(B)$ is dense in
$\op{Pic}^{0}(C)$. Then for every $d \in {\mathbb Z}$ there exist line
bundles ${\mathcal N} \in \op{Pic}^{d}(C)$, s.t. $\T_{C}({\mathcal N})
= {\mathcal N}$.
\end{lem}
{\bf Proof.} The morphism $\T_{C} : \op{Pic}(C) \to \op{Pic}(C)$ is given
explicitly by the formula:
\begin{equation} \label{eq-ni-inv}
\T_{C}({\mathcal N}) = \alpha_{C}^{*}({\mathcal N})\otimes
{\mathcal O}_{C}(e_{9} - e_{1} + f),
\end{equation}
where $\alpha_{C} = \alpha_{B|C}$.

Indeed, by part (b) of \cite[Proposition~7.7]{dopw-i} this formula
holds for all line bundles 
${\mathcal N} \in \op{Pic}(i_{C}^{*}\op{Pic}(B))$. By
the density assumption it holds for all ${\mathcal N} \in
\op{Pic}^{0}(C)$. But applying $\T_{C}$ to the short exact sequence 
\[
0 \to {\mathcal N}(-p) \to {\mathcal N} \to {\mathcal O}_{p} \to 0
\]
we find $\T_{C}({\mathcal N}(-p)) = \T_{C}({\mathcal
N})(-\alpha_{C}(p))$, so the formula extends to all components of
$\op{Pic}(C)$.

Thus a point 
$x \in  \op{Pic}^{0}(C)$ will be fixed under $\T_{C}$ if and
only if
\begin{equation} \label{eq-TCi-fixed}
x - \alpha_{C}^{*}(x) = {\mathcal O}_{C}(e_{9} - e_{1} + f).
\end{equation}
This equation is consistent exactly when 
\[
{\mathcal O}_{C}(e_{9} - e_{1} + f) \in \op{im}[
\xymatrix@1{\op{Pic}^{0}(C)
\ar[rr]^-{\alpha_{C}^{*} - \op{id}} && \op{Pic}^{0}(C)}].
\]
Since $\alpha_{C}$ has fixed points on $C$, it follows
\cite{mumford-prym}  that  
$\op{im}(\alpha_{C}^{*} - \op{id}) = \ker(\alpha_{C}^{*} +
\op{id})$. But from \cite[Table~1~and~Formula~(4.2)]{dopw-i} we see that
\[
(\alpha_{B}^{*} + \op{id})(e_{9} - e_{1} + f) = 2 e_{9} + 2f + e_{7} -
\ell = o_{1} + o_{2}.
\]
Since $o_{1}$ and $o_{2}$ do not intersect $C$, this implies that 
${\mathcal O}_{C}(e_{9} - e_{1} + f)$ is
$\alpha_{C}^{*}$-anti-invariant. Hence there is a translate of
of $\op{Pic}^{0}(C/\alpha_{C})$ in $\op{Pic}^{0}(C)$ consisting of solutions
of \eqref{eq-TCi-fixed}.  

Shifting by an arbitrary multiple of an
$\alpha_{C_{i}}$-fixed point, we see that there are $\T_{C_{i}}$-fixed
points in $\op{Pic}^{d}(C_{i})$ for every $d$. \hfill $\Box$

\

\bigskip

In view of this lemma it only remains to check the density of
$i_{C_{i}}^{*}\op{Pic}(B)$ in  $\op{Pic}^{0}(C_{i})$. We do this only
in the cases $(i = 2, k_{2} = 3)$ and $(i=3, k_{3} = 6)$, which are
the cases needed in section~\ref{s-numerics}. 
Unfortunately the statement of Lemma~\ref{lem-invariant-mod-density}
does not directly apply to the first of these cases (see
Remark~\ref{rem-vertical}), so we will treat it separately next. 
\

\subsection{Invariance for $k_{2} = 3$} \label{sss-C2} By
Remark~\ref{rem-vertical}, the general $\alpha_{B}$-invariant curve
$C_{2}$ in the linear system $|2e + 3f|$ is of the form $C_{2} =
\overline{C}_{2} + F$, where
$\overline{C}_{2}$ is a smooth curve in the linear system $|2e +
2f|^{+}$ and $F$ denotes one of the elliptic curves $f_{0}$,
$f_{\infty}$. 
Assume that
${\mathcal N}_{2} \to C_{2}$ is a line bundle and let
$\overline{{\mathcal N}}_{2} = {\mathcal N}_{2}\otimes {\mathcal
O}_{\overline{C}_{2}}$ be its restriction to $\overline{C}_{2}$. We
know that $\overline{W}_{2} :=
\FM_{B}(i_{\overline{C}_{2}*}\overline{{\mathcal N}}_{2})$ is a vector
bundle. We want $W_{2}:= \FM_{B}(i_{C_{2}*}{\mathcal N}_{2})$ to be a
vector bundle too.

\begin{lem} \label{lem-degree1} $W_{2}$ is a vector bundle if and only
if $\deg({\mathcal N}_{2|F}) = 1$.
\end{lem}
{\bf Proof.} We have a short exact sequence of torsion sheaves on $B$
\[
0 \to i_{F*}({\mathcal N}_{2|F}(-D)) \to i_{C_{2}*}{\mathcal N}_{2}
\to i_{\overline{C}_{2}*}\overline{{\mathcal N}}_{2} \to 0,
\]
where $D \subset F$ is the effective divisor $D =
\overline{C}_{2}\cap F$. Let $G := {\mathcal N}_{2|F}(-D)$. Since $G$
is a line bundle on the fiber $F$ we have $\FM_{B}(i_{F*}G) =
i_{F*}(\FM_{F}(G))$, where $\FM_{F} : D^{b}(F) \to  D^{b}(F)$ is the
Fourier-Mukai transform with respect to the Poincare bundle ${\mathcal
P}_{B|F\times F}$. If we apply $\FM_{B}$ to the
above exact sequence, we will get the 
long exact sequence of cohomology sheaves 
\[
\xymatrix@R=6pt{
0 \ar[r] & {\mathcal H}^{0}(i_{F*}\FM_{F}(G)) \ar[r] & {\mathcal
H}^{0}\FM_{B}(i_{C_{2}*}{\mathcal N}_{2})  \ar[r] & {\overline{W}}_{2}
\ar`r_l/5pt[lll]`^dr/5pt[lll][dll] & \\
&  {\mathcal H}^{1 }(i_{F*}\FM_{F}(G)) \ar[r] & {\mathcal H}^{1
}\FM_{B}(i_{C_{2}*}{\mathcal N}_{2}) \ar[r] 
& 0. &
}
\]
Since we want the line bundle ${\mathcal N}_{2} \to C_{2}$ 
to be chosen so that 
${\mathcal H}^{1}\FM_{B}(i_{C_{2}*}{\mathcal N}_{2}) = 0$ and 
$W_{2} = {\mathcal H}^{0}\FM_{B}(i_{C_{2}*}{\mathcal N}_{2})$ is a
rank two vector bundle on $B$, we must have 
${\mathcal H}^{0}(i_{F*}\FM_{F}(G)) = 0$ and   ${\mathcal
H}^{1}(i_{F*}\FM_{F}(G))$ must be a line bundle on $F$ such that there
exists a surjection $\overline{W}_{2|f} \twoheadrightarrow {\mathcal
H}^{1}(i_{F*}\FM_{F}(G))$. This can only happen if $G$ has  
degree $-1$ on $F$. \hfill $\Box$

\bigskip

\noindent
We note that in the situation of the lemma $W_{2}$ fits in a short
exact sequence 
\[
0 \to W_{2} \to \overline{W}_{2} \to i_{F*}(G^{\vee}) \to 0, 
\]
where $G$ is defined in the proof of the lemma. Indeed, the proof of 
Lemma~\ref{lem-degree1}gave us a short exact sequence 
\[
0 \to W_{2} \to \overline{W}_{2} \to i_{F*}(\FM_{F}(G)[1]) \to 0.
\]
But every line
bundle of degree $-1$ on an elliptic curve $F$ is of the form
${\mathcal O}_{F}(-p)$ for some point $p \in F$. Applying $\FM_{F}$ to
the short exact sequence 
\[
0 \to {\mathcal O}_{F}(-p) \to {\mathcal O}_{F} \to {\mathcal O}_{p}
\to 0
\]
we see that $\FM_{F}(G) = G^{\vee}[-1]$ for any line bundle of degree
$-1$ on $F$.

Given $(\overline{C}_{2},\overline{{\mathcal N}}_{2})$
Lemma~\ref{lem-degree1} produces a 2-parameter family of vector
bundles $W_{2}$. Indeed, let $G \to F$ be any
line bundle of degree $-1$. Consider the semi-stable bundle
$\overline{W}_{2|F}$ on $F$. Since generically $\overline{C}_{2}$
intersects $F$ at two distinct points we will have $\overline{W}_{2|F}
= A\oplus A^{\vee}$, where $A$ is a non-trivial line bundle of degree
zero on $F$. Therefore $h^{0}(F,A^{\vee}\otimes G^{\vee}) =
h^{0}(F,A\otimes G^{\vee}) = 1$ and so we have unique (up to scale)
maps $A \to G^{\vee}$ and $A^{\vee} \to G^{\vee}$. Also since the
degree one line bundles $A^{\vee}\otimes G^{\vee}$ and $A\otimes
G^{\vee}$ are not isomorphic, it follows that their unique (up to a
scale) sections vanish at different points on $F$. Hence we get a one
parameter family of surjective maps of vector bundles $A\oplus
A^{\vee} \to G^{\vee}$. The corresponding Hecke transform of 
$\overline{W}_{2}$:
\[
W_{2} = \ker[\overline{W}_{2} \to i_{F*}(A\oplus
A^{\vee}) \to i_{F*}(G^{\vee})]
\]
is a rank two vector bundle on $B$ which is the Fourier-Mukai image of
a line bundle ${\mathcal N}_{2}$ on $C_{2} = \overline{C}_{2}\cup
F$. In particular $W_{2}$ is a deformation of a rank two vector bundle $W$
corresponding to  a spectral datum $(C,{\mathcal N})$ where $C$ is a
smooth (but non-invariant) 
curve in the linear system $|2e + 3f|$ and ${\mathcal N}$ is a
line bundle on $C$. Even though this $W$ can not be $\tau_{B}$
invariant, this shows that as far as Chern classes are concerned the
bundle $W_{2}$ behaves as a bundle corresponding to a smooth spectral
cover. 

We are now ready to analyze the $\tau_{B}$-invariance properties of
$W_{2}$. First of all, since $\overline{C}_{2}$ is smooth
Lemma~\ref{lem-invariant-mod-density} applies modulo the following
density statement:

\begin{lem} \label{lem-C2-density}
$i_{\overline{C}_{2}}^{*}\op{Pic}(B)$ is Zariski dense in
$\op{Pic}^{0}(\overline{C}_{2})$. 
\end{lem}
{\bf Proof.} The curves 
$\overline{C}_{2}$ have genus $2$ and in the linear system $|2e +
2f|^{+}$ we have a degenerate curve consisting of the zero section $e$
taken with multiplicity two and of two fibers of $\beta$ which are
exchanged by $\alpha_{B}$. In particular the Jacobian of this
degenerate curve is just the product of the two fibers. But the
Mordell-Weil group of $B$ has rank $6$ and in particular we get
elements of infinite order in the general fiber of $\beta$ which are
restrictions of global line bundles. By continuity this implies that
for a general $\overline{C}_{2} \in |2e + 2f|^{+}$ we can find both
$\alpha_{B}$-invariant and $\alpha_{B}$-anti-invariant line bundles on
$B$ that restrict to elements of infinite order in
$\op{Pic}^{0}(\overline{C}_{2})$. Finally $\alpha_{B}$ has two fixed
points on a general $\overline{C}_{2}$ and so
$g(\overline{C}_{2}/\alpha_{\overline{C}_{2}}) = 1$ and $\dim
\op{Prym}(\overline{C}_{2}, \alpha_{\overline{C}_{2}}) = 1$. Hence 
$i_{\overline{C}_{2}}^{*}\op{Pic}(B)$ is dense in 
$\op{Pic}^{0}(\overline{C}_{2})$. \hfill $\Box$

\bigskip

We now reach the main point of this subsection:

\begin{prop} \label{prop-invariantW2} For every integer $d$ there
exists a $\tau_{B}$-invariant vector bundle $W_{2}$ whose spectral
data $(C_{2}, {\mathcal N}_{2} \in \op{Pic}(C_{2}))$ deforms flatly to 
a smooth curve in $B$ and a line bundle of degree $d$ on it. 
\end{prop}
{\bf Proof.} We have a two parameter family of $\alpha_{B}$-invariant
curves $\overline{C}_{2}$.  By Lemmas \ref{lem-invariant-mod-density} and
\ref{lem-C2-density}, there is a one parameter family of
$\T_{\overline{C}_{2}}$-invariant line bundles $\overline{{\mathcal
N}}_{2}$ on each. Althogether we get a three parameter family of
$\tau_{B}$-invariant bundles $\overline{W}_{2}$. We have seen above
that each of these gives rise to a two parameter family of bundles
$W_{2}$. We will check now that in each such two parameter family
there is a finite number (in fact, four) of $\tau_{B}$ invariant $W_{2}$.

Indeed for every such $\overline{W}_{2}$ we must look for a
$\tau_{B}$-invariant Hecke transform $W_{2}$. For this we need to
ensure that $G^{\vee}$ is preserved by $\tau_{B|F}$ and that the map 
$\overline{W}_{2} \to i_{F*}(G^{\vee})$ is $\tau_{B}$-equivariant. We
have two possibilities:  $F = f_{0}$ or $F = f_{\infty}$. Recall that 
$\tau_{B|f_{0}} = t_{\zeta(0)}$ and $\tau_{B|f_{\infty}} =
t_{\zeta(\infty)}\circ (-1)_{f_{\infty}}$ and that $\zeta(0) \in
f_{0}$ and $\zeta(\infty) \in f_{\infty}$ are non-trivial points of
order two. In particular $\tau_{B|f_{0}}$ does not fix any line bundle
of degree one on $f_{0}$ and $\tau_{B|f_{\infty}}$ fixes precisely
four such bundles, namely the four square roots of the degree two
line bundle ${\mathcal O}_{f_{\infty}}(e(\infty) + \zeta(\infty))$.

Choose now $F = f_{\infty}$ and 
$G^{\vee}$ to be one of the four square roots of ${\mathcal
O}_{f_{\infty}}(e(\infty) + \zeta(\infty))$. Choose a 
non-zero map $s : A \to G^{\vee}$. Then $\tau_{B|f_{\infty}}^{*}s :
A^{\vee} \to G^{\vee}$ is also a non-zero map and so, as before,  
$s\oplus \tau_{B|f_{\infty}}^{*}s : A\oplus A^{\vee} \to G^{\vee}$ is
surjective. Using this map as the center of a Hecke transform, we get a 
$\tau_{B}$-invariant $W_{2}$. \hfill $\Box$

\bigskip

\

\subsection{Invariance for $k_{3} = 6$}  \label{sss-C3}
 Let $W_{3} = \FM_{B}(i_{C_{3}*}{\mathcal N}_{3})$ for
some curve in  $|3e + 6f|$. As we saw above, in this case, we can
chose $C_{3}$ to be smooth and preserved by $\alpha_{B}$ and so in
order to find $\tau_{B}$-invariant $W_{3}$'s we only need to show that
for a general $C_{3} \in |3e + 6f|^{\pm}$ the image
$i_{C_{3}}^{*}\op{Pic}(B)$ will be Zariski dense in
$\op{Pic}^{0}(C_{3})$.  

The Zariski closure of the image $i_{C_{3}}^{*}\op{Pic}(B)$ varies
lower-semi-continuously with $C_{3}$, so it suffices to exhibit one
good $C_{3}$. Our $C_{3}$ will be reducible, consiting of a generic
$\alpha_{B}$-invariant curve $\overline{C}_{2}$ in the linear system
$|2e + 2f|$, plus the zero section $e$, plus two generic fibers
$\phi_{1}$ and $\phi_{2}$, plus their images $\phi_{3} :=
\alpha_{B}(\phi_{1})$ and $\phi_{4} := \alpha_{B}(\phi_{2})$. 

The arithmetic genus of $C_{3}$ is $13$. The $13$-dimensional
$\op{Pic}^{0}(C_{3})$ is a $({\mathbb C}^{\times})^{7}$ extension of 
the six dimensional abelian variety ${\mathcal A} :=
\op{Pic}^{0}(\overline{C}_{2})\times \prod_{i =
1}^{4}\op{Pic}^{0}(\phi_{i})$. So our density 
statement follows from the from the following two lemmas.

\begin{lem} \label{lem-densityA} For a generic choice of
$\overline{C}_{2}$, $\phi_{1}$, $\phi_{2}$, the image of 
$i_{C_{3}*}\op{Pic}^{0}(B)$ in ${\mathcal A}$ is Zariski dense.
\end{lem}

\

\begin{lem} \label{lem-density-fibers} For a generic choice of
$\overline{C}_{2}$, $\phi_{1}$, $\phi_{2}$, no proper subgroup of 
$\op{Pic}^{0}(C_{3})$ surjects onto ${\mathcal A}$.
\end{lem}

\

\bigskip

\noindent
{\bf Proof of Lemma~\ref{lem-densityA}.} Under the genericity
assumption in the hypothesis of the lemma, it is clear that 
there are no isogenies
among $\phi_{1}$, $\phi_{2}$, and the two elliptic curves
$\op{Pic}^{0}(\overline{C}_{2}/\alpha_{\overline{C}_{2}})$ and 
$\op{Prym}(\overline{C}_{2},\alpha_{\overline{C}_{2}})$. So it
suffices to prove density separately in each of the two dimensional
abelian varieties $\phi_{1}\times \phi_{3}$, $\phi_{2}\times \phi_{4}$
and $\op{Pic}^{0}(\overline{C}_{2})$. Density in
$\op{Pic}^{0}(\overline{C}_{2})$ was already proved in
Lemma~\ref{lem-C2-density} and density in say $\phi_{1}\times \phi_{3}$
was established during the proof of that result. The lemma is proven.
\ \hfill $\Box$

\

\bigskip

\noindent
{\bf Proof of Lemma~\ref{lem-density-fibers}.} Let ${\mathcal B}$ be a
${\mathbb C}^{\times}$-extension of an abelian variety ${\mathcal
A}$. Such a ${\mathcal B}$ determines a line bundle $[{\mathcal B}]
\in \op{Pic}^{0}({\mathcal A})$. A proper subgroup of ${\mathcal B}$
surjecting onto ${\mathcal A}$ will exist if and only if $[{\mathcal
B}]$ is torsion. 

Similarly, our $({\mathbb C}^{\times})^{7}$-extension
$\op{Pic}^{0}(C_{3})$ will contain a proper subgroup surjecting onto
${\mathcal A}$ if and only there is a non-zero character $\chi :
({\mathbb C}^{\times})^{7} \to {\mathbb C}^{\times}$ such that the
associated line bundle $L_{\chi} := 
\op{Pic}^{0}(C_{3})\times_{\chi}{\mathbb C}$
over ${\mathcal A}$ is torsion. 

Therefore it suffices to find seven characters $\chi_{1}, \ldots
\chi_{7}$ of $({\mathbb
C}^{\times})^{7}$ such that the associated line bundles are linearly 
independent over ${\mathbb Q}$. For this we will need an intrinsic
description of the character lattice $\Lambda$ of $({\mathbb
C}^{\times})^{7}$ in terms of the geometry of the curve $C_{3}$. 
The singular set of $C_{3}$ is $S = \{ s_{ij} | i = 1,2,3,4 \text{ and
} j = 1,2,3 \}$, where 
$\phi_{i}\cap \overline{C}_{2} = \{ s_{ij} \}_{j = 1}^{2}$, and
$s_{i3} =  \phi_{i}\cap e$. Here the singular points are labelled so
that $\alpha_{B}(s_{1j}) = s_{3j}$ and $\alpha_{B}(s_{2j}) = s_{4j}$. 
Now the lattice $\Lambda$ is
explicitly described as:
\[
\Lambda = \ker[{\mathbb Z}^{S} \to \pi_{0}(C_{3} - S)].
\]
Here the map ${\mathbb Z}^{S} \to \pi_{0}(C_{3} - S)$ sends the
characteristic function $\epsilon_{ij}$ of $s_{ij}$ to the difference
$\phi_{i} - 
\overline{C}_{2}$ for $j = 1,2$ and to $\phi_{i} - e$ for $j = 3$. Our 
seven characters $\chi_{k}$ are the following seven elements in 
${\mathbb Z}^{S}$: 
\[
\begin{split}
\chi_{1} & = \epsilon_{11}  + \epsilon_{31} - \epsilon_{12} -
\epsilon_{32}, \\
\chi_{2} & = \epsilon_{21} + \epsilon_{41} - \epsilon_{22} -
\epsilon_{42}, \\
\chi_{3} & = \epsilon_{11} + \epsilon_{33} - \epsilon_{13} -
\epsilon_{31}, \\
\chi_{4} & = \epsilon_{21} + \epsilon_{43} - \epsilon_{23} -
\epsilon_{41}, \\
\chi_{5} & = \epsilon_{11} + \epsilon_{32} - \epsilon_{12} -
\epsilon_{31}, \\
\chi_{6} & = \epsilon_{21} + \epsilon_{42} - \epsilon_{22} -
\epsilon_{41}, \\
\chi_{7} & = \epsilon_{21} + \epsilon_{41} - \epsilon_{11} -
\epsilon_{31} + \epsilon_{13} + \epsilon_{33} - \epsilon_{23} -
\epsilon_{43}.
\end{split}
\]
To prove the independence of the $L_{\chi_{k}}$ over ${\mathbb Q}$ it
suffices to prove the independence of the restrictions
$\lambda_{k} := L_{\chi_{k}|\prod_{i=1}^{4} \phi_{i}}$. 
We represent a degree zero 
line bundle $\lambda$ on $\prod_{i=1}^{4} \phi_{i}$ by a four-tuple 
$(\lambda^{i}
\in \op{Pic}^{0}(\phi_{i}))_{i = 1}^{4}$. In this notation the
$\lambda_{k}$'s become: 
\[
\begin{split}
\lambda_{1} & = (a,0,a,0), \\
\lambda_{2} & = (0,b,0,b), \\
\lambda_{3} & = (c,0,-c,0), \\
\lambda_{4} & = (0,d,0,-d), \\
\lambda_{5} & = (a,0,-a,0), \\
\lambda_{6} & = (0,b,0,-b), \\
\lambda_{7} & = (-c,d,-c,d),
\end{split}
\]
where 
\[
\begin{split}
a & = {\mathcal O}_{\phi_{1}}(s_{11} - s_{12}) \in
\op{Pic}^{0}(\phi_{1}) \cong \op{Pic}^{0}(\phi_{3}), \\
b & = {\mathcal O}_{\phi_{2}}(s_{21} - s_{22}) \in
\op{Pic}^{0}(\phi_{2}) \cong \op{Pic}^{0}(\phi_{4}), \\
c & = {\mathcal O}_{\phi_{1}}(s_{11} - s_{13}) \in
\op{Pic}^{0}(\phi_{1}) \cong \op{Pic}^{0}(\phi_{3}), \\
d & = {\mathcal O}_{\phi_{2}}(s_{21} - s_{23}) \in
\op{Pic}^{0}(\phi_{2}) \cong \op{Pic}^{0}(\phi_{4}). \\
\end{split}
\]
So we only need to show the linear independence over ${\mathbb Q}$ of
$a, c \in \op{Pic}^{0}(\phi_{1})$ and similarly for $b, d \in
\op{Pic}^{0}(\phi_{2})$. This however is obvious from the fact that
$c, d$ are univalued as functions on the base ${\mathbb P}^{1}$,
whereas $a, b$ are two-valued. The lemma is proven. \ \hfill $\Box$

\section{Numerical conditions} \label{s-numerics}

Our goal is to construct a stable rank five holomorphic vector bundle
on the Calabi-Yau manifold $Z := X/\tau_{X}$ which has a trivial
determinant, three generations and an anomaly which can be absorbed
into $M5$-branes. In terms of $X$ this amounts to finding a rank five
vector bundle $V \to X$ so that:

\bigskip

\begin{itemize}
\item[{\bf (S)}] $V$ is a stable vector bundle.
\item[{\bf (I)}] $V$ is $\tau_{X}$-invariant.
\item[{\bf (C1)}] $c_{1}(V) = 0$.
\item[{\bf (C2)}] $c_{2}(X) - c_{2}(V)$ is effective.
\item[{\bf (C3)}] $c_{3}(V) = 12$.
\end{itemize}

\

\bigskip

\noindent
We will construct a whole family of $V$'s satisfying these conditions.
As explained in section~\ref{s-spectral}, each $V$ will be constructed
as an extension
\begin{equation}
0 \to V_{2} \to V \to V_{3} \to 0 \label{eq-extension},
\end{equation}
where the $V_{i}$'s have special form. In fact, as we argued in
section~\ref{s-spectral}, in order to satisfy the condition {\bf (I)}
it is sufficient to take
\[
V_{i} = \FM_{X}\left(\Sigma_{i}, (\pi'_{|\Sigma_{i}})^{*}{\mathcal
N}_{i}\otimes \pi^{*}L_{i} \otimes {\mathcal
O}_{\Sigma_{i}}\left(-\sum_{k=1}^{i}a_{ik}(\{p_{i1k} \} \times n_{1}' +
\{p_{i2k} \} \times o_{2}')\right)\right)
\]
with $a_{ik}$ being positive integers, $\Sigma_{i} = \pi^{'*}C_{i}$
for some smooth curve $C_{i}   \subset B$ satisfying \eqref{eq-ci-inv} and
${\mathcal N}_{i} \in \op{Pic}^{d_{i}}(C_{i})$ satisfying
\eqref{eq-ni-invariance}. In fact, we do not need $C_{i}$ to be
smooth; it is sufficient for the pair $(C_{i},{\mathcal N}_{i})$ to be
deformable to a pair $(C_{i}',{\mathcal N}_{i}')$ with $C_{i}'$ being
smooth but not necessarily $\alpha_{B}$-invariant.
Furthermore, we showed in section~\ref{s-spectral}, that
such $(C_{i},{\mathcal N}_{i})$ do exist and move in positive
dimensional families, at least for specific values of $k_{i}$. 
From now on we will assume that $C_{i}$ and
${\mathcal N}_{i}$  are chosen so that they are deformable to a smooth
pair and that \eqref{eq-ci-inv} and
\eqref{eq-ni-invariance} hold. 

Next we will rewrite the conditions {\bf (S)} and {\bf (C1-3)} as a
sequence of numerical constraints on the numbers $k_{i}$, $d_{i}$ and
on the line bundles $L_{i}$.

\bigskip

\subsection{The Chern classes of $V$} \label{ss-chern-classes}
There are several ways of calculating the Chern classes
of the bundles $V_{i}$. One possibility is to use the cohomological
Fourier-Mukai transform on $X$. To avoid long and cumbersome
calculations of $\fm_{X}$ we choose a slightly different approach
which utilizes the details of the geometric structure of $V$.

Recall that in section~\ref{s-spectral} we gave an alternative
description of the bundles $V_{i}$ as 
\[
V_{i} = \widetilde{W}_{i}\otimes \pi^{*} L_{i},
\]
where $\widetilde{W_{i}}$ is the result of $a_{ik}$, $k = 1, \ldots i$
Hecke transforms of the bundle $\pi^{'*}W_{i}$, where $W_{i} =
\FM_{B}((C_{i},{\mathcal N}_{i}))$. Due to Corollary~\ref{cor-chern}
and the identity \eqref{eq-hecke-pattern} we have
\[
ch(\widetilde{W}_{i}) = \pi^{'*}ch(W_{i}) - \left(\sum_{k =
1}^{i}a_{ik}\right)\pi^{'*}(n_{1}' + o_{2}') - 2\left(\sum_{k =
1}^{i}a_{ik}^{2}\right)(f\times \op{pt}).
\]
Next we need the following
\begin{lem} \label{lem-chernW} Let $C \subset B$ be a smooth curve in
the linear system $|{\mathcal O}_{B}(re + mf)|$ and let ${\mathcal N}
\in \op{Pic}^{d}(C)$. Let $W = \FM_{B}((C,{\mathcal N}))$. Then 
\[
ch(W) = r + \left(d + \begin{pmatrix} r+1 \\ 2 \end{pmatrix} -rm
-r\right)f - m\cdot\op{pt}.
\]
\end{lem}
{\bf Proof.} The bundle $W$ is defined as the Fourier-Mukai transform
of the spectral datum $(C,{\mathcal N})$ on $B$.  Explicitly this
means that $W = \FM_{B}(i_{C*}{\mathcal N})$, where $i_{C} : C
\hookrightarrow B$ is the inclusion map. In particular
\[
ch(W) = \fm_{B}(ch(i_{C*}{\mathcal N}))
\]
and so it suffices to calculate $ch(i_{C*}{\mathcal N})$. 

\medskip

\noindent
By Groethendieck-Riemann-Roch theorem we have
\[
\begin{split}
ch(i_{C*}{\mathcal N}) & = (i_{C*}(ch({\mathcal N})td(C)))td(B)^{-1} \\
& = (i_{C*}(1 + (d + 1 - g)\cdot \op{pt}))\cdot \left(1 + \frac{1}{2} f
+ \op{pt} \right)^{-1} \\
& = ((re + mf) +  (d + 1 - g)\cdot \op{pt}))\cdot \left(1 - \frac{1}{2} f
- \op{pt} \right) \\
& = (re + mf) + \left(d + 1 - g - \frac{r}{2}\right)\cdot \op{pt}.
\end{split}
\]
Also by adjunction we have $2g - 2 =  C\cdot (K_{B} + C) = 2rm - r^{2}
- r$. Hence
\[
ch(i_{C*}{\mathcal N}) = (re + mf) + \left( d + \begin{pmatrix} r+1 \\
2 \end{pmatrix} - rm - \frac{r}{2} \right)\cdot \op{pt}.
\]
Finally using \cite[Table~2]{dopw-i} we get 
\[
\begin{split}
ch(W) & = \fm_{B}(ch(i_{C*}{\mathcal N})) \\
& = \fm_{B}\left((re + mf) + \left( d + \begin{pmatrix} r+1 \\
2 \end{pmatrix} - rm - \frac{r}{2} \right)\cdot \op{pt}\right) \\
& = r \left( 1 - \frac{1}{2}f\right) + m\cdot (-\op{pt}) + \left( d +
\begin{pmatrix} r+1 \\ 
2 \end{pmatrix} - rm - \frac{r}{2} \right)\cdot f \\
& = r + \left( d + \begin{pmatrix} r+1 \\ 2 \end{pmatrix} - rm - r
\right)f - m\cdot \op{pt}.
\end{split}
\]
The lemma is proven. \hfill $\Box$

\

\begin{rem} \label{rem-singularOK} Since the Chern classes are
topological invariants,  the conclusion of the previous lemma still
holds even if $C$ is singular, as long as the sheaf $i_{C*}{\mathcal
N}$ deforms flatly to a line bundle on some smooth spectral curve in
the same linear system. In particular it will hold for the $W_{2}$
from section~\ref{sss-C2}. 
\end{rem}

\bigskip

\noindent
Going back to the calculation of $ch(V)$ let us write $S_{i}^{\alpha}$
for the Newton sums
\[
S_{i}^{\alpha} = \sum_{k=1}^{i} a_{ik}^{\alpha}.
\]
In this notation we need to calculate the product
\[
ch(V_{i}) = (\pi^{'*}ch(W_{i}) - S_{i}^{1}\pi^{*}(n_{1}' + o_{2}') - 2
S_{i}^{2} (f\times \op{pt}))\cdot ch(\pi^{*}L_{i}).
\]
But 
\[
ch(\pi^{*}L_{i}) = \pi^{*}ch(L_{i}) = 1 + \pi^{*}L_{i} +
\frac{L_{i}^{2}}{2}(f\times \op{pt}),
\]
and so  
\begin{itemize}
\item $\pi^{*}(n_{1}' + o_{2}')\cdot ch(\pi^{*}L_{i}) = \pi^{*}(n_{1}'
+ o_{2}') + (L_{i}\cdot (n_{1}' + o_{2}'))\cdot (f\times \op{pt})$;
\item $(f\times \op{pt})\cdot ch(\pi^{*}L_{i}) = f\times \op{pt}$;
\item Lemma~\ref{lem-chernW} yields
\[
\begin{split}
\pi^{'*}& ch(W_{i})\cdot ch(\pi^{*}L_{i}) =\\
& = \pi^{'*}\left(i + 
\left( d_{i} + \begin{pmatrix} i+1 \\ 2 \end{pmatrix} - k_{i}i -
i\right)\cdot f - k_{i}\cdot \op{pt}\right)\cdot \left(1 +
\pi^{*}L_{i} +
\frac{L_{i}^{2}}{2}(f\times \op{pt})\right)\\ 
& = i + \pi^{*}\left( iL_{i} +
\left( d_{i} + \begin{pmatrix} i+1 \\ 2 \end{pmatrix} - k_{i}i -
i\right) f'\right) + \\
& \quad + \left[ \left( \frac{i}{2}L_{i}^{2} + \left( d_{i} +
\begin{pmatrix} i+1 \\ 2 \end{pmatrix} - k_{i}i - 
i\right)\cdot (L_{i}\cdot f')\right)(f\times \op{pt}) -
k_{i}(\op{pt}\times f')\right] \\
& \quad   - k_{i}(L_{i}\cdot f')\op{pt}.
\end{split}
\]
\end{itemize}
Combining these formulas we get formulas for $ch(V_{2})$ and
$ch(V_{3})$:
\[
\boxed{
\begin{split}
ch(V_{2}) & = 2 + \pi^{*}(2L_{2} + (d_{2} - 2k_{2} + 1)f' -
S_{2}^{1}(n_{1}' + o_{2}')) \\
& \quad +  
(L_{2}^{2} + (d_{2} - 2k_{2} + 1)(L_{2}\cdot f') -
S_{2}^{1}(L_{2}\cdot (n_{1}' + o_{2}')) - 2 S_{2}^{2})\cdot (f\times
\op{pt})\\
& \qquad \qquad \qquad - k_{2}(\op{pt}\times f')\\
& \quad -k_{2} (L_{2}\cdot f')\op{pt}.
\end{split}
}
\]
and similarly
\[
\boxed{
\begin{split}
ch(V_{3})  & = 3 + \pi^{*}(3L_{3} + (d_{3} - 3k_{2} + 3)f' -
S_{3}^{1}(n_{1}' + o_{2}')) \\
& \quad +
\left(\frac{3}{2}L_{3}^{2} + (d_{3} - 3k_{3} + 3)(L_{3}\cdot f') -
S_{3}^{1}(L_{3}\cdot (n_{1}' + o_{2}')) - 2 S_{3}^{2}\right)\cdot (f\times
\op{pt}) \\
& \qquad \qquad \qquad - k_{3}(\op{pt}\times f') \\
& \quad -k_{3} (L_{3}\cdot f')\op{pt}.
\end{split}
}
\]
Together these formulas give
\[
\boxed{\boxed{
\begin{split}
ch(V) & = 5 + \\
& + \pi^{*}(2L_{2} + 3L_{3} + (d_{2} + d_{3} - 2k_{2} - 3k_{3} + 4)f' -
(S_{2}^{1} + S_{3}^{1})(n_{1}' + o_{2}')) \\
& + [(L_{2}^{2} + (3/2)L_{3}^{2} + (d_{2} - 2k_{2} + 1)(L_{2}\cdot f')
+ (d_{3} - 3k_{3} + 3)(L_{3}\cdot f') \\
&\qquad \qquad \qquad  - (S_{2}^{1}L_{2} + S_{3}^{1}L_{3})\cdot (n_{1}'
+ o_{2}') - 2(S_{2}^{2} + S_{3}^{2}))(f\times \op{pt}) \\
& \qquad \qquad \qquad \qquad \qquad - (k_{2} + k_{3})(\op{pt}\times f')] \\
& - ((k_{2}L_{2} + k_{3}L_{3})\cdot f')\op{pt}.
\end{split}
}}
\]
Therefore, taking into account that $c_{2}(X) = 12(f\times \op{pt} +
\op{pt}\times f')$ we see that the conditions {\bf (C1-3)} translate
into the following numerical constraints:

\bigskip

\begin{description}
\item[\num{C1}{}] $2L_{2} + 3L_{3} = (S_{2}^{1} + S_{3}^{1})(n_{1}' +
o_{2}') - (d_{2} + d_{3} - 2k_{2} - 3k_{3} + 4)f'$.
\item[\num{C2}{f}] $k_{2} + k_{3} \leq 12$.
\item[\num{C2}{f'}]
$L_{2}^{2} + (3/2)L_{3}^{2} + (d_{2} - 2k_{2} + 1)(L_{2}\cdot f')
+ (d_{3} - 3k_{3} + 3)(L_{3}\cdot f') - (S_{2}^{1}L_{2} + 
S_{3}^{1}L_{3})(n_{1}' + o_{2}') - 2(S_{2}^{2} + S_{3}^{2}) \geq -12$. 
\item[\num{C3}{}] $k_{2}(L_{2}\cdot f') + k_{3}(L_{3}\cdot f') = -6$.
\end{description}

\bigskip

\noindent
Our next task is to express the stability of $V$ in a numerical form.

\subsection{Stability of $V$} \label{ss-stability}

We need to make sure that the bundle $V$ is Mumford stable with
respect to some ample class $H \in H^{2}(X,{\mathbb Z})$. Recall (see
e.g. \cite[Section 7]{fmw-vb}) that
a polarization $H \in H^{2}(X,{\mathbb Z})$ is called {\em suitable}
if up to an overall scale the components of the fibers of $\pi : X \to
B'$ have 
sufficiently small volume. Starting with any polarization $H_{0} \in
H^{2}(X,{\mathbb Z})$ we can construct a suitable polarization by
fixing some polarization $h' \in H^{2}(B',{\mathbb Z})$ and taking 
\[
H := H_{0} + n\cdot \pi^{*}h',
\] 
with $n \gg 0$.  
As explained
in \cite[Theorem~7.1]{fmw-vb} for a suitable $H$ every
vector bundle on $X$ which comes from a spectral cover 
will be $H$-stable on $X$. From now on we will
always assume that $H = H_{0} + n\cdot \pi^{*}h'$ is chosen to be suitable.
For a torsion free sheaf ${\mathcal F}$ on
$X$ denote by  $\mu_{H}({\mathcal F})$ the $H$-slope of ${\mathcal
F}$, i.e.   $\mu_{H}({\mathcal
F}) = (c_{1}({\mathcal F})\cdot H^{2})/\op{rk}({\mathcal F})$.
By repeating the argument in the proof
of \cite[Theorem~7.1]{fmw-vb} one gets the following lemma whose proof
is left to the reader.
\begin{lem} \label{lem-stability}
The bundle $V$ constructed in the previous section is $H$-stable if
and only if
\begin{itemize}
\item[{\em (i)}] The extension 
\[
0 \to V_{2} \to V \to V_{3} \to 0,
\]
is non-split.
\item[{\em (ii)}] $\mu_{H}(V_{2}) < \mu_{H}(V) = 0$. 
\end{itemize}
\end{lem}

\bigskip

\noindent
Next we express both these conditions into a numerical form. Note that
an extension
\[
0 \to V_{2} \to V \to V_{3} \to 0
\]
will be non-split if and only if it corresponds to a non-zero element
in $\op{Ext}^{1}(V_{3},V_{2}) = H^{1}(X,V_{3}^{\vee}\otimes
V_{2})$. Thus we only need to ensure that $H^{1}(X,V_{3}^{\vee}\otimes
V_{2}) \neq 0$. We have the following
\begin{lem} \label{lem-extension} For $V_{2}$ and $V_{3}$ as above one has
$H^{1}(X,V_{3}^{\vee}\otimes V_{2}) \neq 0$ if  $L_{2}\cdot
f' > L_{3}\cdot f'$.
\end{lem}
{\bf Proof.} We are assuming that the $W_{i}$'s deform to vector
bundles on $B$ coming from smooth spectral covers. So by the
upper-semi-continuity of $H^{1}(X,V_{2}^{\vee}\otimes V_{3})$ it is
enough to prove the lemma for $W_{i}$'s arising from smooth spectral
covers.

Let $L = L_{3}^{-1}\otimes L_{2}$. Then 
\[
V_{3}^{\vee}\otimes V_{2} = \widetilde{W}_{3}^{\vee}\otimes
\widetilde{W}_{2}\otimes \pi^{*}L.
\]
To calculate $H^{1}(X,\widetilde{W}_{3}^{\vee}\otimes
\widetilde{W}_{2}\otimes \pi^{*}L)$ we use the Leray spectral
sequence for the projection $\pi : X \to B'$. It 
yields an exact sequence of vector spaces
\[
\xymatrix@C=10pt@R=4pt{
0 \ar[r] &  H^{1}(B',\pi_{*}(\widetilde{W}_{3}^{\vee}\otimes
\widetilde{W}_{2})\otimes L) \ar[r] &  H^{1}(X,V_{3}^{\vee}\otimes
V_{2}) \ar[r] & H^{0}(B',R^{1}\pi_{*}(\widetilde{W}_{3}^{\vee}\otimes
\widetilde{W}_{2})\otimes L) \\ 
\ar[r] & H^{2}(B',\pi_{*}(\widetilde{W}_{3}^{\vee}\otimes
\widetilde{W}_{2})\otimes L). & & 
}
\]
By construction $\widetilde{W}_{2}$ and $\widetilde{W}_{3}$ are vector 
bundles on $X$ coming from the spectral construction applied to line
bundles on {\em smooth} spectral covers, which are also finite over $B'$.
In particular the restriction of $\widetilde{W}_{3}^{\vee}\otimes
\widetilde{W}_{2}$ to the general fiber of $\pi$ is regular and
semistable of degree zero. Hence for general $\widetilde{W}_{2}$ and
$\widetilde{W}_{3}$ we have
\begin{itemize}
\item $\pi_{*}(\widetilde{W}_{3}^{\vee}\otimes \widetilde{W}_{2}) = 0$;
\item $R^{1}\pi_{*}(\widetilde{W}_{3}^{\vee}\otimes
\widetilde{W}_{2})$ is supported on a curve in $B'$ and is a line
bundle on its support.
\end{itemize}
Therefore
\[
H^{1}(X,V_{3}^{\vee}\otimes V_{2}) =
H^{0}(B',R^{1}\pi_{*}(\widetilde{W}_{3}^{\vee}\otimes 
\widetilde{W}_{2})\otimes L).
\]
To calculate the latter space notice that
$R^{1}\pi_{*}(\widetilde{W}_{3}^{\vee}\otimes 
\widetilde{W}_{2})$ is supported on the curve
\[
\pi((-1)_{X}^{*}\Sigma_{3}\cap \Sigma_{2}) \subset B'.
\]
Since $\Sigma_{i} = \pi^{'*}C_{i}$, this implies
\[
(-1)_{X}^{*}\Sigma_{3}\cap \Sigma_{2} = \coprod_{t \in
(-1)_{B}^{*}C_{3}\cap C_{2}} \{ t \}\times f'.
\]
Without a loss of generality we may assume that all $t \in
(-1)_{B}^{*}C_{3}\cap C_{2}$ project to distinct points in $\cp{1}$
under the map $\beta : B \to \cp{1}$.  Consequently
\[
R^{1}\pi_{*}(\widetilde{W}_{3}^{\vee}\otimes 
\widetilde{W}_{2}) = \oplus_{t \in
(-1)_{B}^{*}C_{3}\cap C_{2}} \Phi_{t},
\]
where $\Phi_{t}$ is a line bundle on the curve $\pi(\{ t \}\times f')
= f_{\beta(t)}'$. The line bundle $\Phi_{t}$ depends only on the
the restriction of $\widetilde{W}_{3}^{\vee}\otimes 
\widetilde{W}_{2}$ to the surface $f_{\beta(t)}\times
f_{\beta(t)}'$. Let $\op{pr}_{1} : f_{\beta(t)}\times
f_{\beta(t)}' \to f_{\beta(t)}$ and $\op{pr}_{2} : f_{\beta(t)}\times
f_{\beta(t)}' \to f_{\beta(t)}'$ denote the natural projections. Then
we have
\[
\begin{split}
\Phi_{t} & = R^{1}\op{pr}_{2*}((\widetilde{W}_{3}^{\vee}\otimes 
\widetilde{W}_{2})_{|f_{\beta(t)}\times
f_{\beta(t)}'}) = R^{1}\op{pr}_{2*}(\op{pr}_{1}^{*}(W_{3}^{\vee}\otimes
W_{2})) \\
& = R^{1}\op{pr}_{2*}(\op{pr}_{1}^{*}{\mathcal O}_{f_{\beta(t)}}) =
H^{1}(f_{\beta(t)},{\mathcal O}_{f_{\beta(t)}})\otimes {\mathcal
O}_{f_{\beta(t)}'}.
\end{split}
\]
In other words $H^{1}(X,V_{3}^{\vee}\otimes V_{2}) \neq 0$ if and only
if $L_{|f_{\beta(t)}'}$ is effective for some $t \in 
(-1)_{B}^{*}C_{3}\cap C_{2}$. For this it suffices to have $L\cdot f'
>0$, and it is necessary to have $L\cdot f' \geq 0$. 

However we saw in the previous section
that the condition {\bf (C1)} implies $2L_{2}\cdot f' + 3 L_{3}\cdot
f' = 0$. If we assume that $L\cdot f' = 0$, then we will get
$L_{2}\cdot f' = L_{3}\cdot f' = 0$ which contradicts {\bf (C3)}. Thus
$H^{1}(X,V_{3}^{\vee}\otimes V_{2}) \neq 0$ iff $L\cdot f' > 0$. The
lemma is proven. \hfill $\Box$

\bigskip

\bigskip

\noindent
Expressing the slope condition in numerical terms is completely
straightforward. In the previous section we showed that
\[
c_{1}(V_{2}) = \pi^{*}(2L_{2} + (d_{2} + 1 - 2k_{2})f' -
S_{2}^{1}(n_{1}' + o_{2}')).
\]
Hence $\mu_{H}(V_{2}) < 0$ if and only if 
\[
\pi^{*}(2L_{2} + (d_{2} + 1 - 2k_{2})f' -
S_{2}^{1}(n_{1}' + o_{2}'))\cdot H^{2} < 0.
\]
It is more convenient to rewrite this as a condition on the surface
$B'$. Recall that we take $H$ to be of the form 
\[
H = H_{0} + n\cdot \pi^{*}h'
\]
where $h' \in H^{2}(B',{\mathbb Z})$ is some polarization and $n \gg 0$. 
Since $X = B\times_{\cp{1}} B'$ we see that any polarization $H_{0}$
on $X$ can be written as 
\[
H_{0} = \pi^{'*}h_{0} + \pi^{*}h_{0}'
\]
for some polarizations $h_{0} \in H^{2}(B,{\mathbb Z})$ and $h_{0}'
\in H^{2}(B',{\mathbb Z})$. In particular
\[
\begin{split}
H^{2} & = (\pi^{'*}h_{0} + \pi^{*}h_{0}')^{2} + 2n(\pi^{'*}h_{0} +
\pi^{*}h_{0}')\pi^{*}h' + n^{2}\cdot \pi^{*}(h^{'2}) \\
& = (h_{0}^{2})\cdot (\{\op{pt} \}\times f') + 2(\pi^{*}h_{0}' +
n\cdot \pi^{*}h')\pi^{'*}h_{0} + \pi^{*}((h_{0}' + n\cdot h')^{2}).
\end{split}
\]
By the projection formula we get 
\[
\begin{split}
\mu_{H}(V_{2}) & = \pi^{*}(2L_{2} + (d_{2} + 1 - 2k_{2})f' -
S_{2}^{1}(n_{1}' + o_{2}'))\cdot H^{2} \\
& = (2L_{2} + (d_{2} + 1 - 2k_{2})f' -
S_{2}^{1}(n_{1}' + o_{2}'))\cdot \pi_{*}(H^{2}) \\
& = (2L_{2} + (d_{2} + 1 - 2k_{2})f' -
S_{2}^{1}(n_{1}' + o_{2}'))\cdot ((h_{0}^{2})f' + 2(h_{0}' + 
n\cdot h')\pi_{*}\pi^{'*}h_{0}) \\
& = (2L_{2} + (d_{2} + 1 - 2k_{2})f' -
S_{2}^{1}(n_{1}' + o_{2}'))\cdot ((h_{0}^{2})f' + 2(h_{0}\cdot
f)h_{0}') \\
& \qquad \qquad + 2(h_{0}\cdot f)n(2L_{2} + (d_{2} + 1 - 2k_{2})f' -
S_{2}^{1}(n_{1}' + o_{2}'))\cdot h'. 
\end{split}
\]
To derive the last identity we used \eqref{eq-main-diagram} to write
\[
\pi_{*}\pi^{'*}h_{0} = \beta^{'*}\beta_{*}h_{0} = (h_{0}\cdot f)\cdot
1 + m f' \in H^{\bullet}(B',{\mathbb Z}),
\]
with $m$ being a positive integer. This implies  that
$\alpha\cdot \pi_{*}\pi^{'*}h_{0} =  (h_{0}\cdot f)(\alpha\cdot 1) +
m(\alpha\cdot f')$ for any cohomology class $\alpha$ on $B'$. In
particular $\op{pt}\cdot \pi_{*}\pi^{'*}h_{0} = (h_{0}\cdot f) \op{pt}$
and hence the above formula.

\bigskip

In conclusion, we see that for $n \gg 0$ we have
\[
\mu_{H}(V_{2}) < 2(h_{0}\cdot f)n(2L_{2} + (d_{2} + 1 - 2k_{2})f' -
S_{2}^{1}(n_{1}' + o_{2}'))\cdot h',
\]
and so $\mu_{H}(V_{2}) < 0 $ provided that
\[
(2L_{2} + (d_{2} + 1 - 2k_{2})f' -
S_{2}^{1}(n_{1}' + o_{2}'))\cdot h' < 0.
\]
We are now ready to list all the conditions on $V$ in a numerical form.

\subsection{The list of constraints} \label{ss-list}

In the previous two sections we translated the conditions {\bf (S)},
{\bf (I)}, {\bf (C1-3)} into a set of numerical conditions. Together
those read:

\bigskip

\begin{description}
\item[\num{S}{e}] $L_{2}\cdot f' > L_{3}\cdot f'$.
\item[\num{S}{s}] $(2L_{2} + (d_{2} + 1 - 2k_{2})f' - S_{2}^{1}(n_{1}'
+ o_{2}'))\cdot h' < 0$ for some ample class $h' \in \op{Pic}(B')$.
\item[\num{I}{}] $\tau_{B}^{*}L_{i} = L_{i}$, for $i = 2,3$.
\item[\num{C1}{}] $2L_{2} + 3L_{3} = (S_{2}^{1} + S_{3}^{1})(n_{1}' +
o_{2}') - (d_{2} + d_{3} - 2k_{2} - 3k_{3} + 4)f'$.
\item[\num{C2}{f}] $k_{2} + k_{3} \leq 12$.
\item[\num{C2}{f'}]
$L_{2}^{2} + (3/2)L_{3}^{2} + (d_{2} - 2k_{2} + 1)(L_{2}\cdot f')
+ (d_{3} - 3k_{3} + 3)(L_{3}\cdot f') - (S_{2}^{1}L_{2} + 
S_{3}^{1}L_{3})(n_{1}' + o_{2}') - 2(S_{2}^{2} + S_{3}^{2}) \geq -12$. 
\item[\num{C3}{}] $k_{2}(L_{2}\cdot f') + k_{3}(L_{3}\cdot f') = -6$.
\end{description}

\bigskip

\noindent
Observe that these conditions already constrain severely the possible
values of $k_{2}$, $k_{3}$, $L_{2}\cdot f'$ and $L_{3}\cdot f'$. 
Indeed, intersecting both sides of \num{C1}{} with the curve $f'
\subset B'$ we see that 
$2(L_{2}\cdot f') + 3(L_{3}\cdot f') = 0$. Recall also that we showed
in section~\ref{ss-invariance} that for the existence of smooth curves
$C_{2}$ and $C_{3}$ one needs to take $k_{2} \geq 2$ and $k_{3} \geq
3$. Thus the integers $k_{2}$, $k_{3}$, $L_{2}\cdot f'$ and $L_{3}\cdot
f'$ should satisfy:
\begin{itemize}
\item $L_{2}\cdot f' > L_{3}\cdot f'$
\item $k_{2} \geq 2$ and $k_{3} \geq 3$;
\item $k_{2} + k_{3} \leq 12$;
\item $2(L_{2}\cdot f') + 3(L_{3}\cdot f') = 0$;
\item $k_{2}(L_{2}\cdot f') + k_{3}(L_{3}\cdot f') = -6$.
\end{itemize}
Solving these, we find the following finite list of values for
$k_{2}$, $k_{3}$, $L_{2}\cdot f'$ and $L_{3}\cdot f'$. 

\begin{table}[!ht]
\begin{center}
\begin{tabular}{|c|c||c|c|} \hline
$k_{2}$ & $k_{3}$ & $L_{2}\cdot f'$ & $L_{3}\cdot f'$ \\ \hline\hline
$2$ & $4$ & $9$ & $-6$ \\ \hline
$2$ & $6$ & $3$ & $-2$ \\ \hline
$3$ & $5$ & $18$ & $-12$ \\ \hline
$3$ & $6$ & $6$ & $-4$ \\ \hline
$4$ & $7$ & $9$ & $-6$ \\ \hline
\end{tabular}
\end{center} 
\caption{Possible values of $k_{2}$, $k_{3}$, $L_{2}\cdot f'$,
$L_{3}\cdot f'$}
\label{table-list}
\end{table}

\subsection{Some solutions} \label{ss-solutions}

In this section we show that the numerical constraints $(\#)$ can all
be satisfied. In fact, we find infinitely many solutions of
$(\#)$. These represent an infinite sequence of of moduli spaces (of
arbitrarily large dimension) of all possible $V$'s.

Fix $k_{2}$ and $k_{3}$ from the values in Table~\ref{table-list}. For
such a choice the corresponding numbers $L_{2}\cdot f'$ and
$L_{3}\cdot f'$ are the solutions to the linear system
\[
\begin{pmatrix}
2 & 3 \\
k_{2} & k_{3}
\end{pmatrix} \cdot
\begin{pmatrix}
L_{2}\cdot f' \\ L_{3}\cdot f'
\end{pmatrix} =
\begin{pmatrix}
0 \\ -6
\end{pmatrix}.
\]
Thus in terms of $k_{2}$ and $k_{3}$ we have 
\begin{equation} \label{eq-Ls}
L_{2}\cdot f' = \frac{18}{2k_{3} - 3k_{2}}, \quad 
L_{3}\cdot f' = \frac{-12}{2k_{3} - 3k_{2}}.
\end{equation}
Put $k = 2k_{3} - 3k_{2}$. Express $L_{i}$, $i = 2,3$  in terms of the
standard classes on $B'$ as follows
\[
\begin{split}
L_{2} & = +\frac{9}{k}(e' + \zeta') + x_{2}f' +
y_{2}(n_{1}' + o_{2}') + 3M \\
& \\
L_{3} & = -\frac{6}{k}(e' + \zeta') + x_{3}f' + 
y_{3}(n_{1}' + o_{2}') - 2M.
\end{split}
\]
The most general way to make the $L_{i}$'s satisfy \num{I}{} together
with \eqref{eq-Ls} is to take $M$ to be $\tau_{B}$-invariant and
perpendicular to $e' + \zeta'$, $f'$ and $n_{1}' + o_{2}'$. This
follows from the fact that the intersection form on $H^{2}(B',{\mathbb
Z})$ is non-degenerate on the span of $e' + \zeta'$, $f'$ and $n_{1}'
+ o_{2}'$. This is evident from Table~\ref{table-pairings}

\begin{table}[!ht]
\begin{center}
\begin{tabular}{|c||c|c|c|} \hline
 & $e' + \zeta'$ & $f'$ & $n_{1}' + o_{2}'$ \\ \hline\hline
$e'+\zeta'$ & $-2$ & $2$ & $2$ \\ \hline
$f'$ & $2$ & $0$ & $0$ \\ \hline
$n_{1}' + o_{2}'$ & $2$ & $0$ & $-4$ \\ \hline
\end{tabular}
\end{center} 
\caption{Intersection pairing on $e' + \zeta'$, $f'$ and $n_{1}' + o_{2}'$}
\label{table-pairings}
\end{table}
The condition \num{C1}{} translates into
\begin{equation} \label{eq-d2+d3}
\begin{split}
2x_{2} + 3x_{3} & = - (d_{2} + d_{3} - 2k_{2} - 3k_{3} + 4) \\
2y_{2} + 3y_{3} & = S_{2}^{1} + S_{3}^{1}.
\end{split}
\end{equation}
Using Table~\ref{table-pairings} we compute
\begin{align*}
L_{2}^{2} & = -\frac{162}{k^{2}} - 4y_{2}^{2} + 9M^{2} +
\frac{36}{k}(x_{2} + y_{2}) \\
L_{3}^{2} & = -\frac{72}{k^{2}} - 4y_{3}^{2} + 4M^{2} -
\frac{24}{k}(x_{3} + y_{3}) \\
(n_{1}' + o_{2}')\cdot L_{2} & = +\frac{18}{k} - 4y_{2} \\
(n_{1}' + o_{2}')\cdot L_{3} & = -\frac{12}{k} - 4y_{3}.
\end{align*}
So the condition \num{C2}{f'} becomes
\[
\begin{split}
-\frac{270}{k^{2}} - & 4y_{2}^{2} - 6y_{3}^{2}  + 15M^{2} +
 \frac{36}{k}(x_{2} + y_{2} - x_{3} - y_{3}) + \frac{1}{k}(18d_{2} -
 12d_{3} - 36k_{2} + 36k_{3} -18)  \\
& \\
& - S_{2}^{1}\left( \frac{18}{k} - 4y_{2}\right) - S_{3}^{1}\left(-
 \frac{12}{k} - 4y_{3} \right) - 2(S_{2}^{2} + S_{3}^{2}) \geq -12.
\end{split}
\]
We eliminate $d_{3}$ using \eqref{eq-d2+d3} and complete the squares
involving $y_{2}$ and $y_{3}$, to find:
\[
\begin{split}
-\frac{135}{k^{2}} -4u_{2}^{2} - 6u_{3}^{2} & + 15M^{2} +
 ((S_{2}^{1})^{2} - 2S_{2}^{2}) + 
 \left( \frac{2}{3}(S_{3}^{1})^{2} - 2 S_{3}^{2}\right) \\
& \\
& + 
\frac{30}{k}(2x_{2} + d_{2} - 2k_{2} + 1) \geq -12,
\end{split}
\]
where
\[
\begin{split}
u_{2} & = y_{2} - \frac{1}{2}\left( \frac{9}{k} + S_{2}^{1}\right) \\
& \\
u_{3} & = y_{3} - \frac{1}{3}\left( -\frac{9}{k} + S_{3}^{1}\right).
\end{split}
\]
Implementing the second condition in \eqref{eq-d2+d3} we get $2u_{2} +
3u_{3} = 0$. Introduce new variables
\[
\begin{split}
u & := 2u_{2} = - 3u_{3} \\
x & := 2x_{2} + d_{2} - 2k_{2} + 1.
\end{split}
\]
Substituting back into the expressions for $L_{2}$ and $L_{3}$ we get 
\begin{equation}
\begin{aligned}
   L_2 &= \frac{9}{k}\left(e' + \zeta'\right) 
      + \frac{1}{2}\left(x - d_2 + 2k_2 - 1\right)f'
      + \frac{1}{2}\left(u + \frac{9}{k} + S_2^1\right)
          \left(n_1'+o_2'\right)
      + 3M  \\
   L_3 &= -\frac{6}{k}\left(e' + \zeta'\right) 
      + \frac{1}{3}\left(-x - d_3 + 3k_3 - 3\right)f'
      + \frac{1}{3}\left(- u - \frac{9}{k} + S_3^1\right)
          \left(n_1'+o_2'\right)
      - 2M.
\end{aligned}
\label{eq:Lr}
\end{equation}
Similarly for the conditions \num{C2}{f'} and \num{S}{s} we get 
\begin{equation} \label{eq-c2f'}
\frac{5}{3}u^{2} - 15M^{2}  \leq 12 - \frac{135}{k^{2}} +
\frac{30}{k}x +  ((S_{2}^{1})^{2} - 2S_{2}^{2}) +
 \left( \frac{2}{3}(S_{3}^{1})^{2} - 2 S_{3}^{2}\right),
\end{equation}
and 
\begin{equation} \label{eq-Ss}
\left(\frac{18}{k}(e' +\zeta') + xf' + \left(u +
\frac{9}{k}\right)(n_{1}' + o_{2}') + 6M\right)\cdot h' < 0.
\end{equation}
respectively.

We will use the flexibility we have in choosing $M$ to show that
\eqref{eq-c2f'} and \eqref{eq-Ss} have solutions that lead to {\em
integral} coefficients in \eqref{eq:Lr}. The key observation here is
that since  $\op{span}(e' + \zeta', f', n_{1}' + o_{2}') \subset
H^{2}(B',{\mathbb Z})$ contains an ample class, the Hodge index
theorem implies that $M^{2} \leq 0$. Therefore, one expects that there
will be non-effective admissible $M$'s which will make \eqref{eq-Ss}
easier to satisfy.

Note that the {\em means inequality} implies that $(S_{2}^{1})^{2}
- 2S_{2}^{2} \leq 0$ with equality if and only if all the $a_{2k}$'s
are equal to each other. Similarly $(2/3)(S_{3}^{1})^{2} - 2 S_{3}^{2}
\leq 0$ with equality if and only if all the $a_{3k}$'s are equal to
each other. In particular, for any choice of numbers $u$, $x$, $k_{2}$,
$k_{3}$, $a_{ik}$ which satisfies \eqref{eq-c2f'}, the numbers $u$,
$x$, $k_{2}$, $k_{3}$ will satisfy 
\begin{equation} \label{eq-c2f'-weak}
\frac{5}{3}u^{2} \leq 12 + 15 M^{2} - \frac{135}{k^{2}} +
\frac{30}{k}x
\end{equation}
as well.

\bigskip

\noindent
But from Table~\ref{table-list} we see that $k = 2k_{3} - 3k_{2} > 0$
for all admissible values of $k_{2}$ and $k_{3}$. Combined with the
fact that $e' + \zeta'$ is an effective curve this implies
\[
\left(xf' + \left(u +
\frac{9}{k}\right)(n_{1}' + o_{2}') + 6M\right)\cdot h' < 0.
\]
On the other hand $f' = n_{j}' + o_{j}'$ for $j = 1,2$ and so 
$f'\cdot h' > n_{1}'\cdot h'$ and $f'\cdot h' > o_{1}'\cdot h'$. Thus
it suffices to check that 
\begin{equation} \label{eq-Ss-weak}
\left(\left(x + u + \frac{9}{k}\right)f' + 6M\right)\cdot h' < 0. 
\end{equation}
To make things more concrete recall that the only conditions that we
need to impose on $M$ are that $M$ should be $\tau_{B'}^{*}$-invariant
and that $M$ should be perpendicular to $\op{span}(e'+\zeta', f',
n_{1}' + o_{2}')$. From \cite[Table~1]{dopw-i} we see that the classes
$e_{4}' - e_{5}', e_{4}' - e_{6}', 3\ell' - 2(e_{4}' + e_{5}' + e_{6}')
-3e_{7}'$ constitute a rational basis of the space of such $M$'s. Let
us choose for example $M$ to be of the form
\[
M = z(e_{4}' - e_{5}')
\]
for some integer $z$.
With this choice \eqref{eq-c2f'-weak} and \eqref{eq-Ss-weak} become
\begin{equation} \label{eq-c2f'-weak2}
\frac{5}{3}u^2 + 30z^2 - \frac{30}{k}x + \frac{135}{k^2} - 12  \leq 0, 
\end{equation}
and
\begin{equation} \label{eq-Ss-weak2}
\left(\left(x+u+\frac{9}{k}\right)f' 
        + 6z\left(e'_4-e'_5\right) \right) \cdot h' 
      < 0
\end{equation}
respectively.

Consider next the class $\gamma := (x + u + 9/k)f'
+ 6z(e_{4}' -e_{5}')$. Since the K\"{a}hler cone is dual to the Mori
cone we know that an ample class $h'$ with $\gamma\cdot h' < 0$ will exist
as long as $\gamma$ is not effective. But $\gamma$ satisfies 
$\gamma\cdot e_{4}' = x + u + 9/k - 6z$ and so if $6z > x + u + 9/k$
we will have $\gamma\cdot e_{4}' < 0$. Under this assumption we have
two alternatives: either $\gamma$ is not effective or $\gamma -
e_{4}'$ is effective.  However we have $(\gamma -e_{4}')\cdot f' = -1$
and $f'$ moves, so
$\gamma - e_{4}'$ and hence $\gamma$ can not be effective.

In other words, as a first check for the consistency of the
inequalities \eqref{eq-c2f'} and \eqref{eq-Ss} it suffices to make sure
that in the 3-space with coordinates $(x,u,z)$ one can find points 
between the plane 
\[
6z = x + u + \frac{9}{k}
\]
and the paraboloid
\[
\frac{5}{3}u^2 + 30z^2 - \frac{30}{k}x + \frac{135}{k^2} - 12 = 0.
\] 
If we use the equation of the plane to eliminate $x$ and substitute
the result in \eqref{eq-c2f'-weak2} we obtain the quadratic inequality 
\begin{equation} \label{eq-consistency}
\frac{5}{3}\left(u+\frac{9}{k}\right)^2 
      + 30\left(z-\frac{3}{k}\right)^2 - 12 \leq 0,
\end{equation}
which always has solutions regardless of the value of $k$.

To find an actual solution we will choose a particular value for
$k$. By examining Table~\ref{table-list} we see that the possible
values of $k = 2k_{3} - 3k_{2}$ are $1, 2, 3$ and $6$. Furthermore,
since all the coefficients in \eqref{eq:Lr} must be integers, $k$ has
to divide $\op{gcd}(6,9) = 3$ i.e. we may have either $k = 1$ (which
corresponds to $k_{2} = 3$ and $k_{3} = 5$) or $k = 3$
(which corresponds to $k_{2} = 3$ and $k_{3} = 6$. For concreteness we
choose the second case, i.e.
\[
k_{2} = 3, \qquad k_{3} = 6, \qquad k = 3.
\]
Note that the geometry of this case has already been carefully analyzed in
sections~\ref{sss-C2} and \ref{sss-C3}. In particular we showed that
for these values of $k_{i}$ there are spectral pairs $(C_{i},{\mathcal
N}_{i})$ which are deformable to smooth pairs and which satisfy
\eqref{eq-ci-inv} and \eqref{eq-ni-invariance}. 

To minimize \eqref{eq-consistency} we will take
\[
u = -3, \qquad z = 1, \qquad x = 5,
\]
where the value of $x$ is chosen to satisfy \eqref{eq-Ss-weak2}.

We then have 
\begin{equation}
\begin{aligned}
   L_2 &= 3\left(e' + \zeta'\right) 
      + \frac{1}{2}\left(4 - d_2\right)f'
      + \frac{1}{2}\left(6 + S_2^1\right)
          \left(n_1'+o_2'\right)
      + 3\left(e'_4-e'_5\right)  \\
   L_3 &= - 2\left(e' + \zeta'\right) 
      + \frac{1}{3}\left(16 - d_3\right)f'
      + \frac{1}{3}\left(-6 + S_3^1\right)
          \left(n_1'+o_2'\right)
      - 2\left(e'_4-e'_5\right).
\end{aligned}
\label{eq:Lex}   
\end{equation}
We see that all the coefficients in \eqref{eq:Lex} will be integral as
long as: $d_{2}$ is even, $d_{3}\equiv 1\!\pmod 3$, $S_{2}^{1}$ is even
and  $S_{3}^{1}$ is divisible by $3$.

The inequality \eqref{eq-c2f'} now reads
\begin{equation} \label{eq-c2f'-fate}
-2 \leq \left( (S_2^1)^2 - 2S_2^2 \right) 
         + \left( \frac{2}{3}(S_2^1)^2 - 2S_2^2 \right)
\end{equation}
and the inequality \eqref{eq-Ss} reads
\begin{equation} \label{eq-Ss-fate}
\left(6\left(e'+\zeta'\right) + 5f' + 6\left(e'_4-e'_5\right)\right)\cdot
h' < 0
\end{equation}
Note that \eqref{eq-Ss-fate} does not involve the numbers $a_{ik}$ and
so all the restrictions on the $a_{ik}$'s come from
\eqref{eq-c2f'-fate} and from the integrality of \eqref{eq:Lex}.

In view of the discussion about the means inequality above we see that 
\eqref{eq-c2f'-fate} will be automatically satisfied if we take
all $a_{2k}$'s to be equal to a fixed integer $a_{2} \geq 0$ and all
$a_{3k}$'s  to be equal to another fixed integer $a_{3} \geq 0$. 
Moreover with such a choice we clearly have $S_{2}^{1} = 2a_{2}$ and
$S_{3}^{1} = 3a_{3}$ and so we have infinitely many possibilities for
the numbers $a_{ik}$. Since $d_{2}$ and $d_{3}$ are unconstrained
except for the conditions $d_{2}\equiv0\!\pmod2$ and
$d_{3}\equiv1\!\pmod3$ we see that all the conditions 
{\bf (S)}, {\bf (I)} and {\bf (C1-3)} will be satisfied if we can
prove the following:
\begin{claim} \label{claim-existence-of-h'}
There exists an ample class $h' \in \op{Pic}(B')$ satisfying
\eqref{eq-Ss-fate}.
\end{claim}
{\bf Proof.} As explained above the existence of $h'$ is equivalent to
showing that the class 
\[
6(e'+\zeta') + 5f' + 6(e'_4-e'_5) = 6(e_{1}' +
e_{9}') + 5f' + 6(e_{4}' - e_{5}')
\] 
is not in
the Mori cone of $B'$. First consider the class $\xi' :=  
e_{4}' - e_{5}' + 
e_{9}' + f' \in \op{Pic}(B')$. We have  $\xi^{'2} = -1$ and
$\xi\cdot f' = 1$. So $\xi'$ is an exceptional class on $B'$. It
is well known \cite{lnm777} that on a {\em general} rational elliptic
surface every exceptional class is effective and is a section. Since
our $B'$ is not generic we can't use this statement to conclude that
$\xi'$ is the class of a section. However we have

\begin{lem} \label{lem-xi'} The divisor $\xi'$ satisfies
\[
{\mathcal O}_{B'}(\xi') = c_{1}([e_{4}'] -
[e_{5}']).
\] 
In particular $\xi'$ is effective and is a section of
$\beta' : B' \to \cp{1}$.
\end{lem}
{\bf Proof.}  Let $\underline{\xi}'$ be the section of $B'$ for which 
$[\underline{\xi}'] = [e_{4}'] - [e_{5}']$, that is $c_{1}([e_{4}'] - [e_{5}'])
= {\mathcal O}_{B'}(\underline{\xi}')$. 
Since the group law on the general fiber $f_{t}'$, $t \in
\cp{1}$ is defined in terms of the Abel-Jacobi map and since we have
taken $e_{9}'(t) \in f_{t}'$ to be the neutral element for the group
law it follows that 
\[
{\mathcal O}_{B'}(\underline{\xi}' - e_{9}')_{|f_{t}'} = 
(c_{1}([e_{4}'] - [e_{5}'])\otimes {\mathcal
O}_{B'}(-e_{9}'))_{|f_{t}'}  = {\mathcal O}_{B'}(e_{4}' - e_{5}')_{|f_{t}'} 
\]
for the general $t \in \cp{1}$. Therefore the line bundle 
${\mathcal O}_{B'}(\underline{\xi}' + e_{5}' - e_{4}' - e_{9}')$
must be a combination of vertical divisors on $B'$, i.e. we can write 
\begin{equation} \label{eq-unxi'}
\underline{\xi}' = e_{4}' - e_{5}' + e_{9}'
+ a\cdot f' + b\cdot n_{1} + c\cdot n_{2}
\end{equation}
for some integers $a$, $b$ and $c$. By \cite[Formula~(4.2)]{dopw-i} we have 
$e_{4}'\cdot n_{i}' = e_{5}' \cdot n_{i}' = 0$ and $e_{9}'\cdot n_{i}'
= 1$ for $i = 1,2$. Consider the $I_{2}$ fiber $n_{1}'\cup o_{1}'$ of
$B'$. The smooth part $(n_{1}'\cup o_{1}')^{\sharp} := 
(n_{1}'\cup o_{1}')-(n_{1}'\cap o_{1}')$ of
this fiber is an abelian group isomorphic to ${\mathbb Z}/2\times
{\mathbb C}^{\times}$ with $n_{1}'-(n_{1}'\cap o_{1}')$ being the
connected component of the identity. By definition the section
$\underline{\xi}'$ intersects $n_{1}'\cup o_{1}'$ at a point which is the
difference 
of the points $e_{4}'\cdot (n_{1}'\cup o_{1}') = e_{4}'\cdot o_{1}'$
and $e_{5}'\cdot (n_{1}'\cup o_{1}') = e_{5}'\cdot o_{1}'$ in the
group law of $(n_{1}'\cup o_{1}')^{\sharp}$. Since these two points
belong to the same component of $(n_{1}'\cup o_{1}')^{\sharp}$ and the
group of connected components of $(n_{1}'\cup o_{1}')^{\sharp}$ is
${\mathbb Z}/2$ it follows that $\underline{\xi}'$ intersects $n_{1}'\cup
o_{1}'$ at a point in $n_{1}'$, i.e. $\underline{\xi}'\cdot n_{1}' =
1$. Similarly $\delta'\cdot n_{2}' = 1$.  Therefore, intersecting both
sides of \eqref{eq-unxi'} with $n_{1}'$ and $n_{2}'$ we get
$1 = 1 + b$ and $1 = 1 + c$ respectively. Thus $b = c = 0$. Finally
from the fact that $\underline{\xi}^{'2} = -1$ we compute that $a = 1$ and so
$\underline{\xi}' = \xi'$. The lemma is proven. \hfill $\Box$
\

\medskip

\noindent
In view of the previous lemma we have a section $\xi'$ of $\beta' :
B' \to \cp{1}$ and we need to show that the class
\[
\mu := 6e_{1}' + 6\xi' - f' \in \op{Pic}(B')
\] 
is not in the Mori cone of $B'$. 

Assume that $\mu$ is in the Mori cone.
Note that by the definition of $\xi'$ we have 
$e_{1}'\cdot f' = \xi'\cdot f' = 1$, $e_{1}^{'2} = \xi^{'2} = -1$ and 
$e_{1}'\cdot \xi' = 1$. Then $\mu\cdot \xi' = -1$ and so $\mu - \xi' =
6e_{1}' + 5\xi' - f'$ will also have to be in the Mori cone. But now 
$(\mu - \xi')\cdot e_{1}' = - 2$ and so $\mu - 2e_{1}' - \xi'$ will be
in the Mori cone. Intersecting with $\xi'$ again we get $(\mu -
2e_{1}' - \xi')\cdot \xi' = -2$ and so continuing iteratively we
conclude that $-\xi - f$ must be in the Mori cone which is obviously false. 
This shows that $\mu$ is not in the Mori cone of $B'$ and so $h'$
ought to exist. 

For completeness we will identify an explicit ample class $h'$ on $B'$
with $\mu\cdot h' <0$. We will look for $h'$ of the form
\[
h' = af' + be_{1}' + c\xi'
\]
and will try to adjust the coefficients $a$, $b$ and $c$ so that $h'$
is ample and $\mu\cdot h' < 0$. First we have the following 

\begin{lem} \label{lem-ample}
The  divisor class $h' = af' + be_{1}' + c\xi'$ is ample provided that
$a$, $b$ and $c$ are positive and $a > |b - c|$.
\end{lem}
{\bf Proof.} Assume that $a$, $b$ and $c$ are positive and $a > |b -
c|$. By the Nakai-Moishezon criterion for ampleness
\cite[Theorem~1.10]{hartshorne},  $h'$ will be ample if 
$h^{'2} > 0$ and if $h'\cdot C > 0$ for every irreducible curve $C
\subset B'$. 

Let $C \subset B'$ be an irreducible curve. Then we have two
possibilities: either $C$ is a component of a fiber of $\beta' : B'
\to \cp{1}$, or $\beta' : C \to \cp{1}$ is a finite map. If $\beta' :
C \to \cp{1}$ is finite and $C \neq \xi', e_{1}'$, then $C\cdot f' >
0$, $C\cdot e_{1}' \geq 0$  and $C\cdot \xi' \geq 0$. In particular
the fact that $a$, $b$ and $c$ are positive implies that 
$C\cdot h' > 0$. Hence $h'$ will be ample if we can show that the
intersections  $h^{'2}$, $h'\cdot f'$, $h'\cdot n_{i}'$, $h'\cdot
o_{i}'$, $h'\cdot e_{1}'$ and $h'\cdot\xi'$ are all positive. For this
we calculate 
\[
\begin{split}
h'\cdot e_{1}' & = a + c - b \\
h'\cdot\xi' & = a + b - c \\
h'\cdot f' & = b + c \\
h'\cdot n_{1}' & = h'\cdot n_{2}' = c \\
h'\cdot o_{1}' & = h'\cdot o_{2}' = b \\
h^{'2} & = - b^{2} - c^{2} + 2ab + 2ac + 2bc = b(a + c - b) + c(a + b
- c) + ab + ac,
\end{split}
\]
which are manifestly positive provided that $a > |b - c|$. The lemma
is proven. \hfill $\Box$
\

\medskip

\noindent
Now we see that the condition $\mu\cdot h' < 0$ translates into 
$12a < b + c$, i.e we may take $h' = 25 f' + 144 e_{1}' + 168 \xi'$.
This finishes the proof of the claim. \hfill $\Box$

\section{Summary of the construction} \label{s-summary}

In this section we recapitulate the main points of the
construction. Recall that we want to build a quadruple $(X,H,\tau_{X},V)$
satisfying:
\begin{itemize}
\item[$({\mathbb Z}/2)$] $X$ is a smooth Calabi-Yau 3-fold and
$\tau_{X} : X \to X$ is a freely acting involution. $H$ is a fixed
K\"{a}hler structure (ample line bundle) on $X$
\item[{\bf (S)}] $V$ is an $H$-stable vector bundle of rank five on $X$.
\item[{\bf (I)}] $V$ is $\tau_{X}$-invariant.
\item[{\bf (C1)}] $c_{1}(V) = 0$.
\item[{\bf (C2)}] $c_{2}(X) - c_{2}(V)$ is effective.
\item[{\bf (C3)}] $c_{3}(V) = 12$.
\end{itemize}
\

\bigskip

The construction is carried out in several steps.

\subsection{The construction of $(X,\tau_{X})$} \label{ss-X-and-tauX}

$X$ is built as the fiber product of two rational elliptic surfaces of
special type.

\punkt {\bf Building special rational elliptic surfaces.}
\label{sss-build-res} \quad
Let $\Gamma_{1} \subset \cp{2}$ be a nodal cubic with a node
$A_{8}$. Choose four generic points on $\Gamma_{1}$ and label them
$A_{1}, A_{2}, A_{3}, A_{7}$. Let $\Gamma \subset \cp{2}$ be the
unique smooth cubic which passes trough $A_{1}, A_{2}, A_{3}, A_{7},
A_{8}$ and is tangent to the line $\langle A_{7}A_{i} \rangle$ for $i
= 1, 2, 3$ and $8$. Consider the pencil of cubics spanned by
$\Gamma_{1}$ and $\Gamma$. All cubics in this pencil pass trough 
$A_{1}, A_{2}, A_{3}, A_{7}, A_{8}$ and are tangent to $\Gamma$ at
$A_{8}$. Let $A_{4}, A_{5}, A_{6}$ be the remaining three base points,
and let $B$ denote the blow-up of $\cp{2}$ at the points $A_{i}$, $i =
1,2, \ldots, 8$ and the point $A_{9}$ which is infinitesimally near
$A_{8}$ and corresponds to the line $\langle A_{7}A_{8} \rangle$.

The pencil becomes the anti-canonical map $\beta : B \to \cp{1}$ which
is an elliptic fibration with a section. The map $\beta$ has two
reducible fibers $f_{i} = n_{i}\cup o_{i}$, $i = 1,2$ of type
$I_{2}$. We denote by $e_{i}$, $i = 1, \ldots, 7$ and $e_{9}$ the
exceptional divisors corresponding to $A_{i}$, $i = 1, \ldots, 7$ and
$A_{9}$, and by $e_{8}$ the reducible divisor $e_{9} + n_{1}$. The
divisors $e_{i}$ together with the pullback $\ell$ of a class of a line
from $\cp{2}$ form a standard basis  in $H^{2}(B,{\mathbb Z})$.

The surface $B$ has an involution $\alpha_{B}$ which is uniquely
characterized by the properties: $\alpha_{B}$ commutes with $\beta$, 
$\alpha_{B}$ induces an involution on $\cp{1}$, and $\alpha_{B}$ fixes
the proper transform of $\Gamma$ pointwise. 

Choosing  $e_{9}$ as the zero section of $\beta$, we can interpret
any other section $\xi$ as an automorphism $t_{\xi} : B \to B$
which acts along the fibers of $\beta$. The automorphism $\tau_{B} =
t_{e_{1}}\circ \alpha_{B}$ is again an involution of $B$ which
commutes with $\beta$, induces the same involution on $\cp{1}$ as
$\alpha_{B}$ and has four isolated fixed points sitting on the same
fiber of $\beta$.

The special rational elliptic surfaces form a four dimensional
irreducible family. Their geometry was the subject of \cite{dopw-i}.

\punkt {\bf Building $(X,\tau_{X})$.} \label{sss-build-X} \quad
Choose two special rational elliptic surfaces $\beta : B \to \cp{1}$
and $\beta' : B' \to \cp{1}$ so that the discriminant loci of $\beta$
and $\beta'$ in $\cp{1}$ are disjoint, $\alpha_{B}$ and $\alpha_{B'}$
induce the same involution on $\cp{1}$ and the fix loci of $\tau_{B}$
and $\tau_{B'}$ sit over different points in $\cp{1}$. The fiber
product $X := B\times_{\cp{1}} B'$ is a smooth Calabi-Yau  3-fold
which is elliptic and has a freely acting involution $\tau_{B}\times
\tau_{B'}$ and another (non-free) involution $\alpha_{X} :=
\alpha_{B}\times \alpha_{B'}$. For concreteness we fix the elliptic
fibration of $X$ to 
be the projection $\pi : X \to B'$ to $B'$. 

The Calabi-Yau's form a nine dimensional irreducible family.

\punkt {\bf Building $H$.} \label{sss-build-H} \quad Choose any ample
divisor $H_{0}$ on $X$ and take $H = H_{0} + n\cdot \pi^{*}h'$ for
some positive integer $n$. Then the divisor $H$ will be
ample as long as $h'$ is ample on $B'$ and $n \gg 0$. 

Choose $h' = 25 f' + 144 e_{1}' + 168 \xi'$ with $\xi'$ being the
unique section of $\beta' : \cp{1} \to B'$ satisfying $[e_{4}'] -
[e_{5}]' = [\xi']$. The divisor class $h' \in \op{Pic}(B')$ is ample
on $B'$ by Lemma~\ref{lem-ample}.

\subsection{The construction of $V$} \label{ss-V} 

The bundle $V$ is build as a non-split extension 
\[
0 \to V_{2} \to V \to V_{3} \to 0
\]
of two $\tau_{X}$-invariant stable vector bundles $V_{2}$ and $V_{3}$
of ranks $2$ and $3$ respectively.

Each $V_{i}$ is constructed via the spectral cover construction on
$X$.

\punkt {\bf Building $V_{2}$ and $V_{3}$.} \label{sss-V2-and-V3} \quad
Choose curves $\overline{C}_{2}, C_{3} \subset B$, so that
\begin{itemize}
\item $\overline{C}_{2} \in |{\mathcal O}_{B}(2e_{9} + 2f)|$ and $C_{3} \in
|{\mathcal O}_{B}(3e_{9} + 6 f)|$. 
\item $\overline{C}_{2}$, $C_{3}$ are $\alpha_{B}$-invariant.
\item $\overline{C}_{2}$ and $C_{3}$ are smooth and irreducible.
\end{itemize}
Set $C_{2} :=
\overline{C}_{2} + f_{\infty}$ where $f_{\infty}$ is the smooth fiber
of $\beta$ containing the fixed points of $\tau_{B}$.

The space of such $\overline{C}_{2}$'s is an open set (see
section~\ref{ss-invariance}) in 
${\mathbb P}(H^{0}(B,{\mathcal O}_{B}(2e_{9} +
2f))^{+})$ where $H^{0}(B,{\mathcal O}_{B}(2e_{9} +
2f))^{\pm}$ denote the spaces of invariants/anti-invariants for the
$\alpha_{B}$ action on $H^{0}(B,{\mathcal O}_{B}(2e_{9} + 2f)$. Using
the explicit equations \eqref{eq-C2andC3} of the spectral curves  we easily
that all such $C_{2}$ form a $2$ dimensional irreducible family.
The space of permissible $C_{3}$'s is an open subset in the disjoint
union of the projective spaces ${\mathbb P}(H^{0}(B,{\mathcal
O}_{B}(3e_{9} + 6f))^{\pm})$, which have dimensions $8$ and $6$.

Fix an even integer $d_{2}$ and an integer $d_{3}$ satisfying
$d_{3}\equiv 1\!\pmod 3$. Choose line bundles ${\mathcal N}_{2}$ and
${\mathcal N}_{3}$ on $C_{2}$ and $C_{3}$ respectively, which satisfy
\begin{itemize}
\item ${\mathcal N}_{i} \in \op{Pic}^{d_{i}}(C_{i})$ for $i = 1,2$,
\item ${\mathcal N}_{i} = \T_{C_{i}}({\mathcal N}_{i}) :=
\alpha_{B|C_{i}}^{*}{\mathcal N}_{i}\otimes {\mathcal O}_{C_{i}}(e_{1}
- e_{9} + f)$ for $i = 1, 2$.  
\end{itemize}
As explained in sections~\ref{sss-C2} and \ref{sss-C3}, 
such ${\mathcal N}_{2}$'s and ${\mathcal N}_{3}$'s 
are parameterized by  abelian subvarieties of
$\op{Pic}^{d_{2}-1}(\overline{C}_{2})$ and 
$\op{Pic}^{d_{3}}(C_{3})$ of dimensions equal to the genera of the
quotient curves $\overline{C}_{2}/\alpha_{\overline{C}_{2}}$ and 
$C_{3}/\alpha_{C_{3}}$ respectively. Thus there is a one
dimensional space of ${\mathcal N}_{2}$'s and a six dimensional space
of ${\mathcal N}_{3}$'s. 

Let $\Sigma_{i} = C_{i}\times_{\cp{1}} B'$ for $i = 1,2$. Recall that
$\beta' : B' \to \cp{1}$ has two $I_{2}$ fibers $f_{1}'$ and
$f_{2}'$. Let $F_{1}$, $F_{2}$ be the corresponding (smooth) fibers of $\beta :
B \to \cp{1}$. Let $C_{i}\cap F_{j} = \{ p_{ijk} \}_{k = 1}^{i}$ for
$i = 2,3$, $j = 1, 2$. Then $\Sigma_{i} \to C_{i}$ is an elliptic
surface having $2i$ fibers of type $I_{2}$: $\{ p_{ijk}\}\times
f_{j}'$. Also $\Sigma_{i} \subset X$ and the natural projection
$\pi_{|\Sigma_{i}}  : \Sigma_{i} \to B'$ is finite of degree $i$.

Fix non-negative integers $a_{2}$ and $a_{3}$. Define
\[
V_{i} = \FM_{X}\left(\left( \Sigma_{i},
(\pi_{|\Sigma_{i}})^{*}{\mathcal N}_{i}\otimes {\mathcal
O}_{\Sigma_{i}}\left(-a_{i}\sum_{k = 1}^{i} (\{p_{i1k}\}\times n_{1}'
+ \{p_{i2k}\}\times o_{2}')\right)\right)\right) \otimes \pi^{*}L_{i},
\]
where $L_{2}$ and $L_{3}$ are the line bundles
\[
\begin{split}
L_{2} & = 3(e_{1}' + e_{4}' - e_{5}' + e_{9}') + \frac{1}{4}(4 -
d_{2})f' + (3 - a_{2})(n_{1}' + o_{2}') \\
L_{3} & = -2(e_{1}' + e_{4}' - e_{5}' + e_{9}') + \frac{1}{3}(16 -
d_{3})f' + (-2 + a_{3})(n_{1}' + o_{2}'),
\end{split}
\]
on $B'$.

\punkt {\bf Building $V$.} \label{sss-build-V} \quad  Take $V$ to be a
non-split extension of $V_{2}$ by $V_{3}$ which is $\tau_{X}$-invariant. 
As explained in section~\ref{ss-stability}, the space of all such
extensions of $V_{2}$ by $V_{3}$ is the union of projective spaces 
${\mathbb P}(H^{1}(X,V_{3}^{\vee}\otimes V_{2})^{+})\cup 
{\mathbb P}(H^{1}(X,V_{3}^{\vee}\otimes V_{2})^{-})$ where
$H^{1}(X,V_{3}^{\vee}\otimes V_{2})^{\pm}$ denote the
invariants/anti-invariants for the $\tau_{X}$ action on 
$H^{1}(X,V_{3}^{\vee}\otimes V_{2})$. 

Furthermore, it is shown in section~\ref{ss-stability} that
\[
\dim  H^{1}(X,V_{3}^{\vee}\otimes V_{2}) \geq 
(((-1)_{B}^{*}C_{3})\cdot C_{2})\cdot (L_{2}'\cdot f' - L_{3}'\cdot f') =
150,
\]
and from the explicit description of $H^{1}(X,V_{3}^{\vee}\otimes
V_{2})$ we see that the $\pm$ decomposition  breaks this as $150 = 70
+ 80$, so  the dimension of the admissible extensions of $V_{2}$ by
$V_{3}$ is at least $79$.

In other words, for a fixed $(X,\tau_{X},H)$ as above we find
infinitely many components of the moduli space of $V$'s satisfying
({\bf S}), ({\bf I}) and ({\bf C1-3}). Each component corresponds to
a  choice of the integers $a_{2}, a_{3}, d_{2}$ and $d_{3}$
and has dimension $2 +  8 + 1 + 6 + 79 = 96$. 

\

\bigskip
\bigskip

\appendix

\Appendix{Hecke transforms} \label{app-hecke}

\

\bigskip

\noindent
In this appendix we review the definition and some basic properties of the
Hecke transforms (aka `elementary modifications') 
of vector bundles along divisors. For more details the reader may wish
to consult \cite{maruyama1,maruyama2}, \cite{friedman-vb-book}.

\subsection{Definition and basic properties} \label{as1}

Let $X$ be a smooth complex projective variety. Let $\imath : D 
\hookrightarrow X$ be a divisor with normal crossings.

Let $E \rightarrow X$ be a vector bundle and let $(\xi )$ be a short exact
sequence of vector bundles on $D$ of the form
\[ (\xi ) \;\;\; : \;\;\; 0 \longrightarrow F \longrightarrow E_{|D} 
\longrightarrow G \longrightarrow 0. \]
There are two Hecke transforms $\heck^{\pm}_{(\xi)}(E)$  attached to
the pair $(E,(\xi ))$. 

\begin{defi} \label{d1}
\ 

\noindent
\begin{list}{{{\bf(\roman{rom})}}}{\usecounter{rom}}
\item 
The down-Hecke transform of $E$ along $(\xi )$ is the coherent sheaf
$\hdown{\xi}{E}$ de\-fi\-ned by the exact sequence
\[ 0 \longrightarrow \hdown{\xi}{E} \longrightarrow E \longrightarrow
\imath_{*}G \longrightarrow 0. \] 
\item
The up-Hecke transform of $E$ along $(\xi )$ is the coherent sheaf 
\[ \hup{\xi}{E} = (\hdown{\xi^{\vee}}{E^{\vee}})^{\vee}\]
\end{list}
\end{defi}

\medskip

The first properties of the Hecke transforms are given by the following 
two lemmas.

\medskip 

\begin{lem} \label{l1}
\ 

\noindent
\begin{list}{{{\em\bf(\alph{rom})}}}{\usecounter{rom}}
\item The sheaves $\hup{\xi}{E}$ and $\hdown{\xi}{E}$ are locally free.
\label{l1part1}
\item The up-Hecke transform $\hup{\xi}{E}$ of $E$ along $(\xi )$ fits in the
exact sequence
\[ 0 \longrightarrow E \longrightarrow \hup{\xi}{E} \longrightarrow 
\imath_{*}F\otimes 
{\mathcal O}_{X}(D)\longrightarrow 0. \] \label{l1part2}
\item $\hup{\xi}{E}_{|D}$ and $\hdown{\xi}{E}_{|D}$ are furnished with natural exact 
sequences
\[ (\xi^-) \;\;\; : \;\;\; 0 \longrightarrow G(-D) \longrightarrow 
\hdown{\xi}{E}_{|D} \longrightarrow F \longrightarrow 0, \]
\[ (\xi^+) \;\;\; : \;\;\; 0 \longrightarrow G \longrightarrow 
\hup{\xi}{E}_{|D} \longrightarrow F(D) \longrightarrow 0. \]
\label{l1part3}
\item $\hup{\bullet}{\bullet }$ and $\hdown{\bullet}{\bullet}$ are mutually 
inverse in the sense that
\[ \hup{\xi^-}{\hdown{\xi}{E}} = E, \;\; \hdown{\xi^+}{\hup{\xi}{E}} = E. \]
\item $\hdown{\xi}{E} = \hup{\xi}{E}(-D)$.
\end{list}
\end{lem}
{\bf Proof.} The proof of {\bf (a)} is straightforward. For the proof of {\bf 
(b)} recall that by definition the dual bundle $\hup{\xi}{E}^{\vee}$ fits 
in the short exact sequence of sheaves
\[ 0 \longrightarrow \hup{\xi}{E}^{\vee} \longrightarrow E^{\vee} 
\longrightarrow \imath_{*}(F^{\vee}) \longrightarrow 0. \]
Application of $\homsh{{\mathcal O}_{X}}{\bullet}{{\mathcal O}_{X}}$ combined 
with the fact that $\imath_{*}(F^{\vee})$ is torsion yields the 
long exact sequence

\[
\begin{array}{lclclcl}
0 & \rightarrow & \homsh{{\mathcal O}_{X}}{E^{\vee}}{{\mathcal O}_{X}} & \rightarrow
& \homsh{{\mathcal O}_{X}}{\hup{\xi}{E}^{\vee}}{{\mathcal O}_{X}}& \rightarrow & \\
& \rightarrow & \extsh{1}{{\mathcal O}_{X}}{\imath_{*}(F^{\vee})}{{\mathcal O}_{X}} & 
\rightarrow
& \extsh{1}{{\mathcal O}_{X}}{E^{\vee}}{{\mathcal O}_{X}} & 
\rightarrow & \ldots
\end{array}
\]
Furthermore 
$E^{\vee}$ and ${\mathcal O}_{X}$ are both locally free and hence every 
extension of $E^{\vee}$ by ${\mathcal O}_{X}$ splits locally yielding 
$\extsh{1}{{\mathcal O}_{X}}{E^{\vee}}{{\mathcal O}_{X}} = 0$. Thus we obtain the exact
sequence
\[ 0 \longrightarrow E \longrightarrow \hup{\xi}{E} \longrightarrow 
\extsh{1}{{\mathcal O}_{X}}{\imath_{*}(F^{\vee})}{{\mathcal O}_{X}} \longrightarrow
0. \]
To calculate $\extsh{1}{{\mathcal O}_{X}}{\imath_{*}(F^{\vee})}{{\mathcal O}_{X}}$ 
consider the ideal sequence of the divisor $D$:
\[ 0 \longrightarrow {\mathcal O}_{X} \longrightarrow {\mathcal O}_{X}(D) 
\longrightarrow {\mathcal O}_{D}(D) \longrightarrow 0.\]
After apllying $\homsh{{\mathcal O}_{X}}{\imath_{*}(F^{\vee})}{\bullet}$ and 
taking into account that ${\mathcal O}_{X}$ and ${\mathcal O}_{X}(D)$ are locally
free sheaves, we get the exact sequence

\[
\begin{array}{lclclcl}
0 & \longrightarrow & 
\homsh{{\mathcal O}_{X}}{\imath_{*}(F^{\vee})}{{\mathcal O}_{D}{D}} &
\longrightarrow & \extsh{1}{{\mathcal O}_{X}}{\imath_{*}(F^{\vee})}{{\mathcal O}_{X}}
& \longrightarrow & \\
& \longrightarrow &
  \extsh{1}{{\mathcal O}_{X}}{\imath_{*}(F^{\vee})}{{\mathcal O}_{X}(D)} &
\longrightarrow & \ldots
\end{array} 
\] 
To understand the map 
\begin{equation} \label{eq1}
 \extsh{1}{{\mathcal O}_{X}}{\imath_{*}(F^{\vee})}{{\mathcal O}_{X}} 
\longrightarrow \extsh{1}{{\mathcal O}_{X}}{\imath_{*}(F^{\vee})}{{\mathcal O}_{X}(D)},
\end{equation}
consider a point $p \in D \subset X$. Let $R := {\mathcal O}_{X,p}$ and let
$t \in R$ be a local equation of $D$ around $p$. Let $M$ be the finitely
generated $R/tR$ module whose sheafification gives $F^{\vee} \rightarrow D$
in a neighborhood of $p$.

An element $(\alpha) \in \extsh{1}{{\mathcal O}_{X}}{\imath_{*}(F^{\vee})}{{
\mathcal O}_{X}}_{p}$ of the stalk of 
$\extsh{1}{{\mathcal O}_{X}}{\imath_{*}(F^{\vee})}{{\mathcal O}_{X}}$ at $p$ is an
extension of $R$-modules of the form
\[ (\alpha)\;\;\; : \;\;\; 0 \longrightarrow R \longrightarrow A 
 \longrightarrow M \longrightarrow 0,\]
where $M$ is given its $R$-module structure via $R \rightarrow R/tR$. 

The image $(\beta) \in \extsh{1}{{\mathcal O}_{X}}{\imath_{*}(F^{\vee})}{{
\mathcal
 O}_{X}(D)}_{p}$ of $(\alpha)$ under the map (\ref{eq1}) is just the pushout 
of the extension $(\alpha)$ via the homomorphism
\[ R \longrightarrow \frac{1}{t}R. \]
That is, there is a commutative diagram
\[
\xymatrix{
(\alpha) & 0 \ar[r] & R \ar[r] \ar[d] & A 
\ar[r]^{\pi} \ar[d] & M \ar[r] \ar@{=}[d] & 0 \\
(\beta) & 0 \ar[r] & \frac{1}{t}R \ar[r]  & B
\ar[r] & M \ar[r] & 0}
\]  
and $B = (A\oplus\frac{1}{t}R)/R$.

On the other hand, since $tR$ annihilates $M$ we have $\pi(tx) = t\pi(x) 
= 0$, for all $x \in A$. In particular the map 
\[
\begin{array}{llcl}
s : & A\oplus \frac{1}{t}R & \longrightarrow & \frac{1}{t}R \\
 & & & \\
    & x\oplus \frac{f}{t} & \longrightarrow & \frac{tx}{t} + 
\frac{f}{t},
\end{array}
\]
is well defined and descends to $B$ to a map splitting the exact 
sequence
\[ 0 \longrightarrow \frac{1}{t}R \longrightarrow B \longrightarrow M 
\longrightarrow 0.
\]
Therefore the map
\[ \extsh{1}{{\mathcal O}_{X}}{\imath_{*}(F^{\vee})}{{\mathcal O}_{X}} 
\longrightarrow 
\extsh{1}{{\mathcal O}_{X}}{\imath_{*}(F^{\vee})}{{\mathcal O}_{X}(D)}
\]
is the zero map and we get an isomorphism
\[
\xymatrix{
\homsh{{\mathcal O}_{X}}{\imath_{*}(F^{\vee})}{\imath_{*}{\mathcal O}_{D}
\otimes {\mathcal O}_{X}(D)} \ar[r]^\simeq \ar@{=}[d] &
\extsh{1}{{\mathcal O}_{X}}{\imath_{*}(F^{\vee})}{{\mathcal O}_{X}} \\
\imath_{*}F\otimes {\mathcal O}_{X}(D) &}
\]
which concludes the proof of the lemma.
\hfill $\Box$

\bigskip

\bigskip

There is a natural symmetry between the up and down Hecke transforms.
If $X$, $D$, $E$ and $(\xi)$ are as above, then we can form the dual
exact sequence

\[(\xi^{\vee}) \;\;\; : \;\;\;  0 \longrightarrow G^{\vee} \longrightarrow 
E^{\vee}_{|D} \longrightarrow F^{\vee} \longrightarrow 0, \] 
and the up and down Hecke transforms of $E^{\vee}$ along $(\xi^{\vee})$.
The relation with the Hecke transforms of $E$ is given by the following
lemma.

\begin{lem} \label{l2}
\[
\begin{array}{lcl}
\hup{\xi}{E}^{\vee} &  \simeq &  \hdown{\xi^{\vee}}{E^{\vee}}  \\
 & & \\
\hdown{\xi}{E}^{\vee} & \simeq & \hup{\xi^{\vee}}{E^{\vee}}
\end{array}
\]
\end{lem}
{\bf Proof.} Clear.
\hfill $\Box$

\bigskip

\bigskip

\bigskip

\subsection{Geometric interpretation - flips} \label{as2}

Let $(E,\xi)$ be as in Section \ref{as1} and let $\tau \rightarrow 
\p{E}$ be the relatively ample tautological line bundle. Denote
by $Y := \bl{\p{F}}{\p{E}}$ the blow-up of $\p{E}$ along $\p{F}$.

Let $p : Y \rightarrow \p{E}$ be the blow-up morphism and let 
${\mathcal E} \subset Y$ be the exceptional divisor. The image of $Y$ 
under the full linear system $p^{*}\tau\otimes 
{\mathcal O}_{Y}(-{\mathcal E})$ is again a projective bundle $\p{E'} 
\rightarrow X$. We have the following diagram
\[
\xymatrix{
\tau  \ar@{-}[dr] & & Y \ar[dl]_p 
\ar[dr]^{p'} & & \tau' \ar@{-}[dl] \\
& \p{E} \ar[dr]_{f}  & & \p{E'} \ar[dl]^{f'} \\
 & & X & & }
\]
where $\tau \rightarrow \p{E}$ and $\tau' \rightarrow \p{E'}$ are 
relatively ample line bundles having the properties
\[
\begin{array}{lcl}
 f_{*}\tau &  = & E^{\vee}  \\
 & & \\
 f'_{*}\tau' & = & {E'}^{\vee} \\
 & & \\
 {p'}^{*}\tau' & = & p^{*}\tau\otimes {\mathcal O}_{Y}(- {\mathcal E})
\end{array}
\]
To identify $E'$ in terms of Hecke transforms consider the ideal 
sequence of ${\mathcal E}$:
\[ 0 \longrightarrow {\mathcal O}_{Y}(- {\mathcal E}) \longrightarrow 
{\mathcal O}_{Y} \longrightarrow {\mathcal O}_{\mathcal E} \longrightarrow 0. \]
Tensoring by $p^{*}\tau$ we get 
\begin{equation} \label{eq2}
 0 \longrightarrow {p'}^{*}\tau' \longrightarrow p^{*}\tau 
\longrightarrow  p^{*}\tau\otimes {\mathcal O}_{{\mathcal E}} \longrightarrow 0.
\end{equation}
Let $\pi : Y \rightarrow X$ be the composition $\pi = f\circ p = f'\circ p'$.
Consider the $\pi$ direct image of (\ref{eq2}):
\[
0 \longrightarrow \pi_{*}{p'}^{*}\tau' \longrightarrow \pi_{*}p^{*}\tau 
\longrightarrow \pi_{*}(p^{*}\tau\otimes {\mathcal O}_{\mathcal E}) 
\longrightarrow R^{1}\pi_{*}{p'}^{*}\tau' \longrightarrow \ldots 
\]
Observe first that every fiber of $\pi$ is either a projective space or has 
two irreducible components (meeting transversally) each of which is a 
projective space. Furthermore ${p'}^{*}\tau'$ restricted on a component $P$ of 
the fiber is either ${\mathcal O}_{P}(1)$ or ${\mathcal O}_{P}$ and hence by Serre's
vanishing theorem doesn't have higher cohomology. Thus by the base change and
cohomology theorem $R^{1}\pi_{*}{p'}^{*}\tau' = 0$. Next 

\[ \pi_{*}{p'}^{*}\tau' = f'_{*}{p'}_{*}{p'}^{*}\tau' = 
f'_{*}\tau' = {E'}^{\vee}. \]
Here we used that $p' : Y \longrightarrow \p{E'}$ has connected fibers.

Similarly $\pi_{*}p^{*}\tau = E^{\vee}$ and we get
\[ 
0 \longrightarrow {E'}^{\vee} \longrightarrow E^{\vee} \longrightarrow 
\pi_{*}(p^{*}\tau\otimes {\mathcal O}_{\mathcal E}) \longrightarrow 0.
\]
But $\pi_{*}(p^{*}\tau\otimes {\mathcal O}_{\mathcal E}) = \imath_{*}\bar{f}_{*}
(\tau_{|\p{F}})$ where $\bar{f} : \p{F} \longrightarrow D$ is the natural
projection. Hence we get the short exact sequence
\[ 0 \longrightarrow {E'}^{\vee} \longrightarrow E^{\vee} \longrightarrow
\imath_{*}F^{\vee} \longrightarrow 0 \]
and thus
\[  E' = \hup{\xi}{E}. \]
\

\bigskip

\bigskip

\bigskip

\subsection{An example} \label{as3}

Let $X$ and $D$ be as before. One has the short exact sequence:

\[
\xymatrix{
(D) \;\; : & 0 \ar[r]  & N^{\vee}_{D}X \ar@{=}[d]
\ar[r] & \Omega^{1}_{X|D} \ar[r] & \imath_{*}\Omega^{1}_{D} 
\ar[r] & 0 \\
& & {\mathcal O}_{D}(-D) & & &}
\]
We can form the up and down Hecke transforms of $\Omega^{1}_{X}$ along $(D)$.

\begin{lem} \label{l3}
Denote by $\Omega^{1}_{X}(\log D)$ the sheaf of one forms on $X$ with 
logarithmic poles along $D$. Then the up and down Hecke transforms of
$\Omega^{1}_{X}$ along $(D)$ can be identified as follows
\[
\begin{array}{lcl}
\hup{D}{\Omega^{1}_{X}} & = & \Omega^{1}_{X}(\log D) \\
& & \\
\hdown{D}{\Omega^{1}_{X}} & = & \Omega^{1}_{X}(\log D)\otimes 
{\mathcal O}_{X}(-D).
\end{array}
\]
\end{lem}
{\bf Proof.}
To prove the first equality observe that $\Omega^{1}_{X}(\log D)$ fits in
the residue sequence
\[ 0 \longrightarrow \Omega^{1}_{X} \longrightarrow \Omega^{1}_{X}(\log D)
\longrightarrow \imath_{*}{\mathcal O}_{D} \longrightarrow 0, \]
where the map $\Omega^{1}_{X}(\log D) \rightarrow \imath_{*}{\mathcal O}_{D}$
is given by the residue along $D$. On the other hand, according to 
Lemma~\ref{l1} we have an exact sequence
\[ 0 \longrightarrow \Omega^{1}_{X} \longrightarrow \hup{D}{\Omega^{1}_{X}}
\longrightarrow \imath_{*}{\mathcal O}_{D} \longrightarrow 0 \]
and it is easy to check that the two extension classes coincide. \hfill $\Box$

\

\bigskip
\bigskip


\end{document}